\documentclass[11pt,onecolumn]{article}
\usepackage{amscd}
\usepackage{times}
\usepackage{amsmath}
\usepackage{amssymb}
\usepackage{xspace}
\usepackage{theorem}
\usepackage{graphicx}
\usepackage{ifpdf}
\usepackage{url,hyperref}
\usepackage{latexsym}
\usepackage{euscript}
\usepackage{xspace}
\usepackage[all]{xy}
\usepackage{color}
\usepackage{makeidx}
\usepackage{picins,wrapfig}

\long\def\remove#1{}

\newtheorem{theorem}{Theorem}[section] 

\newtheorem{obs}[theorem]{Observation}
\newtheorem{corollary}[theorem]{Corollary}

\newtheorem{definition}[theorem]{Definition}
\newtheorem{proposition}[theorem]{Proposition}

\newcommand {\mm}[1] {\ifmmode{#1}\else{\mbox{\(#1\)}}\fi}

\newcommand{\coker} {\mathrm {coker}}

\newcommand{\img} {\mathrm {img}}

\newcommand{\XX}		{{X}}

\newcommand{\RR} {\rho}
\newcommand{\JJ} {{\mathcal J}}
\newcommand{\BB} {{\mathcal B}}

\newcommand{\cancel}[1]

\begin{document}

\title{Topological Persistence for Circle Valued Maps}

\author{
Dan Burghelea\thanks{
Department of Mathematics,
The Ohio State University, Columbus, OH 43210,USA.
Email: {\tt burghele@math.ohio-state.edu}}
\quad\quad
Tamal K. Dey\thanks{
Department of Computer Science and Engineering,
The Ohio State University, Columbus, OH 43210, USA.
Email: {\tt tamaldey@cse.ohio-state.edu}}
}
\date{}
\maketitle
\begin{abstract}
We study {\it circle valued maps} and consider the 
{\it persistence  of the homology of their fibers}.
The outcome is a finite collection of computable invariants 
which answer the basic questions 
on persistence and in addition encode the topology of the source space 
and its relevant subspaces. Unlike persistence of real valued
maps, circle valued maps enjoy a different class of invariants
called {\em Jordan cells} in addition to bar codes.
We establish a relation between the homology 
of the source space and of its relevant subspaces with these invariants
and provide a new algorithm to compute these invariants from an input
matrix that encodes a circle valued map on an input simplicial complex. 

\vskip .2in

\end{abstract}
\thispagestyle{empty}
\setcounter{page}{0}
\newpage

\section {Introduction} \label{I}
Data analysis provides plenty of scenarios where one ends up 
with a nice space, most often a simplicial complex, a 
smooth manifold, or a stratified space
equipped with a real valued or a circle valued map.
The persistence theory, introduced in~\cite{ELZ02},
provides a great tool for 
analyzing real valued maps with the help of homology. Similar theory for circle valued maps
has not yet been developed in the literature. 
The work in~\cite{SJ} brings the concept of circle valued maps
in the context of persistence by 
deriving a circle valued map for a given data
using the existing persistence theory.
In contrast, we develop a persistence theory
for circle valued maps.

One place where circle valued
maps appear naturally is the area of dynamics of vector fields.
Many dynamics are described by vector fields which admit a minimizing action  (in mathematical terms a Lyapunov closed one form). Such actions can be interpreted as $1$- cocycles 
which are
intimately connected to circle valued maps as shown
in~\cite{BD10}. Consequently, a notion
of persistence for circle valued maps also provides
a notion of persistence for $1$-cocycles which appear 
in some data analysis problems~\cite{JLY,YY}.
In summary,
persistence theory for circle valued maps promises to
play the role for some vector fields as does 
the standard persistence theory for the scalar 
fields~\cite{CEH07,CEH09,ELZ02,ZC05}. 


One of the main concepts of the persistence theory is the
notion of  {\em bar codes}~\cite{ZC05}--invariants 
that characterize a real valued map at the homology level.
The angle (circle) valued maps, when characterized at 
homology level, require 
a new  invariant called {\em Jordan cells} in addition to 
the refinement of the bar codes into four types.

The standard persistence~\cite{ELZ02,ZC05} which we refer
as {\em sublevel persistence} deals with the change in the homology of the sublevel sets which can not make sense for a circle valued map.
However, the change in the homology of the level sets can be considered for   both real  and circle valued maps. The notion of persistence, when considered for the level sets of a real valued map~\cite{DW07}  is referred here as  {\em level persistence}. It refines the sublevel persistence. 
The zigzag persistence introduced  in~\cite{CSD09} provides  complete
invariants (bar codes) for level persistence of (tame) real valued maps. They are defined  using 
representation theory for linear quivers. 

The change in homology of the level sets of a (tame) circle valued map is more complicated because of the 
{\em return } of the level to itself when one goes along the circle.
It turns out that representation theory of cyclic quivers provides the
complete invariants  for persistence in the homology of the level sets of the circle valued maps. This notion of persistence is  called here  the
{\em persistence for circle valued maps} and its invariants, 
{\em bar codes} and {\em Jordan cells} are
shown to be effectively computable.


Our results include a derivation of the homology for the source
space and its relevant subspaces in terms of the invariants 
(Theorem~\ref{T1} and~\ref{T2}). The result also applies
to real valued maps as they are special 
cases of the circle valued maps. This leads to a result
(Corollary~\ref{C1}) which to our knowledge has not yet
appeared in the literature \footnote {it was brought to our attention 
by David Cohen-Steiner that the  extended persistence proposed
in~\cite{CEH09} allows similar connections between homology of source
spaces and persistence.}. A number of other topological results 
which can not be derived from any of the previously defined persistence theories  are described in 
\cite{BH12} providing additional motivation for this work.

After developing the results on invariants, we propose a new algorithm
to compute the bar codes and Jordan cells. 
For a simplicial complex, the entire
computation can be done by manipulating the original matrix that
encodes the input complex and the map. 
The algorithm first builds a block matrix from the original incidence matrix
which encodes linear maps induced in homology among regular and critical level sets, more precisely the quiver representations $\rho_r$ described in section \ref {S4}. 
Next, it iteratively reduces this new matrix eliminating and hence computing
the bar codes. The resulting matrix which is invertible can be further processed to Jordan 
canonical form~\cite{DS58} providing Jordan cells.
The algorithm
for zigzag persistence~\cite{CSD09} when applied to what we refer in  
section 3 as  the {\it infinite cyclic  covering map $\tilde f$} can
compute bar codes but not Jordan cells.
In contrast, our method can 
compute the bar codes and Jordan cells 
simultaneously by manipulating matrices and can also be used as an 
alternative to compute the bar codes in zig-zag persistence. 

\paragraph{Notations.}
We list here some of the notations that are used throughout.
\begin{itemize}
\item For $r$th homology group of a topological space
$X$ under an a priori fixed field $\kappa$, we write 
$H_r(X)$ instead of $H_r(X;\kappa)$.
 
\item For a map $f:X\to Y$ and $K\subseteq Y$ we write $X_K:= f^{-1}(K)$.
 
\item We use $\mathbb Z_{\geq 0}$ and $\mathbb Z_{>0}$
for non-negative and positive integers respectively.
 
\item In our exposition, we need to use open, semi-open, and closed
intervals denoted as $(a,b)$, $(a,b]$ or $[a,b)$, and
$[a,b]$ respectively. To denote an interval,
in general, we use the notation $\{a,b\}$ where "$\{$" 
stands for either ``$[$" or ``$($".
 
\item For a linear map $\alpha: V\to W$ between two vector spaces we write :
\[ 
\ker \alpha:=\{v\in V \mid \alpha(v)=0\},~~
\img\  \alpha:= \{w\in \alpha(V)\subseteq W\},~~
\coker\ \alpha:= W/\alpha(V).
\]

\item A matrix $A$ is said to be  in {\em column echelon} form if
all zero columns, if any, are on the right to nonzero ones 
and the leading entry (the first nonzero number from below) 
of a nonzero column is always strictly below of the leading  
entry of the next column.
Similarly, $A$ is said to be in {\em row echelon} form  if 
all zero rows, if any, are below nonzero ones and the leading  
entry (the first nonzero number from the right) of a nonzero row 
is always strictly to the right of the leading entry of the row below it.

 
If $A$ is an $m\times n$ matrix ($m$ rows and $n$ columns), there exist
an invertible $n\times n$ matrix $R(A)$ 
and an invertible $m\times m$ matrix $L(A)$ so that
$A\cdot R(A)$ is in column echelon form and
$L(A)\cdot A$ is in row echelon form.
Algorithms for deriving the column and row echelon form 
can be found in standard books on linear algebra.

\end{itemize}

\section{Definitions and background}\label{S1}
We begin with the technical definition of tameness of a map.

For  a continuous map $f: X\to Y$ between two topological spaces
$X$ and $Y$, let $X_U= f^{-1}(U)$ for $U\subseteq Y$.
When $U=y$ is a single point, the set $X_y$ is called
a {\em fiber} over $y$ and is also commonly known as the level set of $y$.
We call the  continuous map $f:X\to Y$ {\em good} if 
every $y\in Y$ has a contractible neighborhood $U$ so that the 
inclusion $X_y\to  X_U$ is a homotopy equivalence.
The  continuous map $f:X\to  Y$ is a {\em fibration} if each 
$y\in Y$ has a  neighborhood $U$ so that the maps $f: X_U \to U$ 
and $pr: X_y\times U\to U$ are  fiber wise homotopy equivalent. 
This means that there exist continuous maps $l:X_U\to  X_y\times U$ with 
$pr|_U \cdot l|_U = f|_U$
which, when restricted to the  fiber  for any $z\in U$, 
are  homotopy equivalences. In particular, $f$ is good.

\begin{definition}
A proper continuous map $f:X\to  Y$ is {\em tame} if it 
is good, and for some discrete closed subset  $S\subset Y$,
the restriction $f: X\setminus  f^{-1}(S) \to Y \setminus S$ is a fibration.
The points in $S\subset Y$ which prevent $f$ to be a fibration  
are called {\em critical values}. 
\end{definition}
  
 If $Y= \mathbb R$  and $X$ is compact or $Y= {\mathbb S}^1,$ \    \footnote{ since the map $f$  is proper and $\mathbb S^1$ compact, so is $X$ } then
the set of critical values is finite, say  $s_1 < s_2 < \cdots s_k.$ 
The fibers above them, $X_{s_i},$ are referred to as  {\em singular fibers}. 
All other fibers are called {\em regular}.
In the case of ${\mathbb S}^1$, $s_i$ can be taken as angles and we can 
assume that $0< s_i \leq 2\pi.$  Clearly, for the open interval
$(s_{i-1},s_i)$ the map $f: f^{-1}(s_{i-1}, s_{i}) \to (s_{i-1}, s_{i})$ 
is a fibration  which  implies that  
all fibers over angles in $(s_{i-1}, s_i)$ 
are homotopy equivalent with a fixed regular fiber, say $X_{t_i}$,
with $t_i\in (s_{i-1}, s_i)$.  

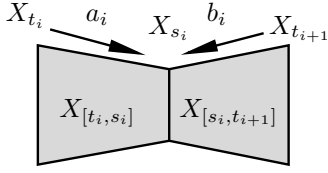
\parpic[l]{
\begin{picture}(0,0)%
\includegraphics{maps.pstex}%
\end{picture}%
\setlength{\unitlength}{3947sp}%
\begingroup\makeatletter\ifx\SetFigFont\undefined%
\gdef\SetFigFont#1#2#3#4#5{%
  \reset@font\fontsize{#1}{#2pt}%
  \fontfamily{#3}\fontseries{#4}\fontshape{#5}%
  \selectfont}%
\fi\endgroup%
\begin{picture}(2532,1054)(1695,-1058)
\put(2026,-136){\makebox(0,0)[b]{\smash{{\SetFigFont{10}{12.0}{\rmdefault}{\mddefault}{\updefault}{\color[rgb]{0,0,0}$X_{t_i}$}%
}}}}
\put(2476,-136){\makebox(0,0)[b]{\smash{{\SetFigFont{10}{12.0}{\rmdefault}{\mddefault}{\updefault}{\color[rgb]{0,0,0}$a_i$}%
}}}}
\put(3751,-211){\makebox(0,0)[b]{\smash{{\SetFigFont{10}{12.0}{\rmdefault}{\mddefault}{\updefault}{\color[rgb]{0,0,0}$X_{t_{i+1}}$}%
}}}}
\put(3226,-136){\makebox(0,0)[b]{\smash{{\SetFigFont{10}{12.0}{\rmdefault}{\mddefault}{\updefault}{\color[rgb]{0,0,0}$b_i$}%
}}}}
\put(2926,-211){\makebox(0,0)[b]{\smash{{\SetFigFont{10}{12.0}{\rmdefault}{\mddefault}{\updefault}{\color[rgb]{0,0,0}$X_{s_i}$}%
}}}}
\put(3301,-736){\makebox(0,0)[b]{\smash{{\SetFigFont{10}{12.0}{\rmdefault}{\mddefault}{\updefault}{\color[rgb]{0,0,0}$X_{[s_i,t_{i+1}]}$}%
}}}}
\put(2476,-736){\makebox(0,0)[b]{\smash{{\SetFigFont{10}{12.0}{\rmdefault}{\mddefault}{\updefault}{\color[rgb]{0,0,0}$X_{[t_i,s_i]}$}%
}}}}
\end{picture}%

}

In particular, there exist maps $a_i: X_{t_i}\to X_{s_i}$ 
and $b_i: X_{t_{i+1}}\to X_{s_{i}}$, unique up to homotopy
defined as follows: If $t_i$ and $t_{i+1}$ are contained in 
$U_i \subset Y$  where the inclusion
$X_{s_i} \subset X_{U_i}$ is a homotopy equivalence with a homotopy inverse $r_i: X_{U_i} \to X_{s_i},$ then $a_i$ and $b_i$ are the restrictions of  $r_i$ to $X_{t_i}$ and $X_{t_{i+1}}$ respectively. If not, in view of the tameness of $f,$ one can find $t'_i$ and $ t'_{i+1}$ in $U_i$  so that $X_{t_i}$ and $X_{t_{i+1}}$ are homotopy  equivalent to $X_{t'_i}$ and $X_{t'_{i+1}}$ respectively
and compose the restrictions of $r_i$  with these homotopy equivalences.  
These maps  determine homotopically  $f:X\to Y,$ when $Y= \mathbb R$ or 
${\mathbb S^1}$. 
For simplicity in writing, when $Y= \mathbb R$ 
we put $t_{k+1}\in (s_k, \infty)$ and $t_1\in (-\infty, s_1)$ and 
when  $Y={\mathbb S}^1$ we put  $t_{k+1}=t_1 \in (s_k, s_1+2\pi).$
All scalar or circle valued simplicial maps 
on a simplicial complex, and all smooth maps with generic 
isolated critical points on a smooth manifold or stratified space
are tame. In particular, Morse maps are tame.

\subsection{Persistence and invariants for real valued maps}\label{S2}
Since our goal is to extend
the notion of persistence from real valued maps 
to circle valued maps, we first summarize the questions 
that the persistence answers when applied to
real valued maps, and then develop a notion
of persistence for circle valued maps which can 
answer similar questions and more. 
We fix a field $\kappa$ and write $H_r(X)$ to denote
the homology vector space of $X$ in dimension $r$ 
with coefficients in a field $\kappa.$

\paragraph{Sublevel persistence.} 
The persistent homology introduced in~\cite{ELZ02} and further developed 
in~\cite{ZC05} is concerned with the following questions:
\begin{itemize}
\item[Q1.]  Does the  class 
$x\in H_r(X_{(-\infty,t]})$ originate in 
$H_r(X_{(-\infty,t'']})$ for $t''< t$?
Does the  class 
$x\in H_r(X_{(-\infty,t]})$ vanish in $H_r(X_{(-\infty,t']})$ for  $t<t'$?
\item[Q2.]  What are the smallest $t'$ and largest $t''$ 
such that this happens?
\end{itemize}

This information is contained in the inclusion induced linear maps 
$H_r(X_{(-\infty,t]}) \to H_r(X_{(-\infty,t']})$ where $t'\geq t$ and 
is known as persistence. Since the involved subspaces are 
sublevel sets, we refer
to this persistence as {\it sublevel persistence}. 
When $f$ is tame, the persistence
for each $r= 0, 1,\cdots \dim X,$ 
is  determined by a finite collection of 
invariants  referred to as {\it bar codes}~\cite{ZC05}.
For sublevel persistence  the bar codes are a collection of 
{\em closed intervals} of the form $[s,s'] $ or 
$[s,\infty)$ with $s, s'$ being the critical values of $f.$ 
From these bar codes one can derive the Betti numbers of $X_{(-\infty, a]},$ 
the dimension of  $\rm{img}(H_r(X_{(-\infty,t]}) \to H_r(X_{(-\infty,t']}))$ 
and get the answers to questions Q1 and Q2.  
For example, the number of $r$-bar codes which contain the interval  
$[a,b]$ is the dimension  of  
$\rm{img}(H_r(X_{(-\infty,a]}) \to H_r(X_{(-\infty,b]})).$ 
The number of $r$-bar codes which identify   to the interval $[a,b]$ is the 
maximal number of linearly independent homology classes born exactly 
in $X_{(-\infty,a]}$ but not before 
and die exactly in $H_r(X_{-\infty, b]})$ but not before. 

\paragraph{Level persistence.}
Instead of sublevels, if we use levels, we obtain what we call
level persistence. The level persistence 
was first considered in~\cite{DW07} but was better understood computationally 
when the zigzag persistence was introduced in~\cite{CSD09}. 
Level persistence is concerned with the homology of the 
fibers $H_r(X_t)$  and addresses questions of the following type. 

\begin{itemize}
\item[Q1.] Does the image of $x\in H_r(X_t)$ vanish in 
$H_r(X_{[t,t']}),$ where $t'>t$ or in $H_r (X_{[t'', t]}),$ where $t''<t$?
\item[Q2.] Can $x$ be detected in $H_r(X_{t'})$ where $t'>t$ 
or in $H_r(X_{t''})$ where $ t'' <t$? 
The precise meaning of  {\it detection} is explained below.
\item[Q3.] What are the smallest $t'$ and the largest  $t''$ 
for the answers to Q1 and Q2 to be affirmative?
\end{itemize}

To answer such questions  one needs information
about the following inclusion induced linear maps:
$$H_r(X_t)\rightarrow H_r(X_{[t, t']})\leftarrow H_r(X_{t'}).$$ 
The {\it level persistence} is the information provided by this collection of  
vector spaces and linear maps for all $t,t'.$ 

We say that $x\in H_r(X_t)$ is dead in $H_r(X_{[t,t']}), \ t'>t ,$
if its image by
$H_r(X_t)\rightarrow H_r(X_{[t, t']})$ vanishes. Similarly, $x$ is dead 
in $H_r(X_{[t'',t]}), \ t''<t ,$ if its image 
by  $H_r(X_t)\rightarrow H_r(X _{[t'', t]})$ vanishes.

\begin{figure}[h]
\centerline{\input{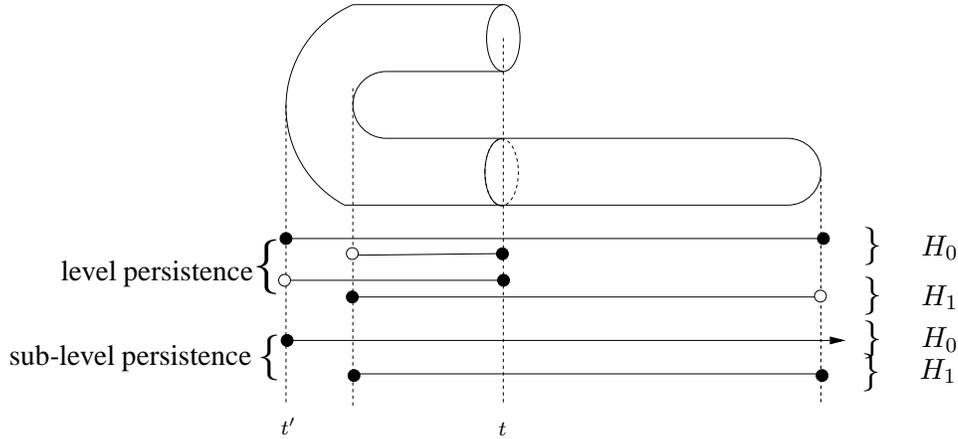}}
\caption{Bar codes for level and sub-level persistence.}
\label{barcode}
\end{figure}

We say that $x\in H_r(X_t)$ is detected   in $H_r (X_{t'}),\  t'>t$, 
 (resp. $t''<t$), if 
its image in $H_r(X_{[t,t']})$ (resp. in $H_r(X_{[t'',t]}$) 
is nonzero and is contained in the image of
$H_r(X_{t'})\rightarrow H_r(X_{[t, t']})$ 
(resp. $H_r(X_{t"})\rightarrow H_r(X_{[t", t]})$). 
In Figure~\ref{barcode}, the class consisting of the sum of two circles
at level $t$ is not detected on the right, but is detected
at all levels on the left up to (but not including) the level $t'$.
In case of a tame map  the collection of the vector spaces 
and linear maps  is 
determined up to coherent isomorphisms by a collection of 
invariants called {\em bar codes for level persistence} which  are intervals of the 
form $[s, s'], (s,s'), (s,s'], [s,s')$ with $s, s'$  
critical values as opposed to the {\em bar codes for sublevel persistence} which are intervals of the form ${[s,s'], [s,\infty)}$ with $s,s' $ critical values.  These bar codes  are called
{\it invariants} because two tame maps $f:X\to \mathbb R$ and $g:Y\to \mathbb R$ which are fiber wise homotopy equivalent have the same associated bar codes.  In the case of level persistence the open end of an interval signifies 
the death of a homology class at that end (left or right) whereas
a closed end signifies that a homology class 
cannot be detected beyond this level (left or right). In the case of the sublevel persistence the left end signifies {\em birth} while the right {\em death}.
Level persistence  provides  considerably more information 
than the sub level persistence. The bar codes of the sub level persistence 
 can 
be recovered from the ones of level persistence. Precisely a  level bar code $[s,s']$ gives a sublevel bar code $[s,\infty)$  and a level bar code $[s,s')$
gives a  sublevel bar code $[s, s'];$  the sublevel persistence does not see any of the level bar codes $(s, s')$ or $(s,s']$.
It turns out that the bar codes of the level persistence can 
also be recovered from the bar codes of the sub level persistence  of $f$  and  additional maps canonically associated to $f.$ 

In Figure~\ref{barcode}, we indicate the bar codes both for 
sub level and level persistence \footnote {the white circles indicate open ends
and the dark circles indicate closed ends} for some simple map $f:X \to \mathbb R$ in order to illustrate their differences. The space $X$ is a tube
open on one end and $f$ is the
height function laid horizontally.

\section{Persistence for circle valued maps}\label{S3}
Let $f:X\to {\mathbb S}^1$ be a circle valued map.
The sublevel persistence for such a map
cannot be defined since circularity in values prevents defining sub-levels.
Even level persistence cannot be defined as per se since
the intervals may repeat over values.
To overcome this difficulty we associate the infinite cyclic covering map  
$\tilde f:\tilde X\to \mathbb R$ for $f$.
It is defined by the commutative diagram: 
\[
\begin{CD}
\tilde \XX                           @>\tilde f>> \mathbb R\\
@V\psi VV                                                       @Vp VV\\
\XX @>f >>            \mathbb S^1
\end{CD} 
\]
The map $p: \mathbb R\to \mathbb S^1$ is the universal 
covering of the circle 
(the map which assigns to the number 
$t\in \mathbb R$ the angle $\theta= t (mod~2\pi)$
and $\psi$ is the pull back of $p$ by the map 
$f$ which is an infinite cyclic covering. 
Notice that if $p(t)=\theta$ then $ \tilde X_t $ and $X_\theta$ are identified by $\psi.$ 
If  $x\in H_r(X_\theta)= H_r(\tilde X_t) ,\  p(t)=\theta,$ 
the questions Q1, Q2, Q3 for $f$ and $X$
can be formulated in terms of the level persistence for $\tilde f$
and $\tilde X$. 

Suppose that $x\in H_r(\tilde X_t)= H_r(X_\theta) $
is detected in $H_r(\tilde X_{t'})$ for some  $t'\geq t+2\pi$.
Then, in some sense, $x$ returns to
$H_r(X_\theta)$ going along the circle ${\mathbb S}^1$
one or more times. When this happens, the class $x$
may change in some respect . This gives rise
to new questions that were not encountered in sublevel
or level persistence.
\begin{itemize}
\item[Q4.] When $x\in H_r(X_\theta)$ returns, how does the
``returned class" compare with the original class $x$?~
It may disappear after going  along the circle a number of times, 
or it might never disappear and if so how does this class change
after its return.
\end{itemize} 
To answer Q1-Q4 one has to record information about
$
H_r(X_\theta)\rightarrow H_r(X_{[\theta, \theta']})\leftarrow H_r(X_{\theta'})
$
for any pair of angles $\theta$ and $\theta'$ 
which differ by at most $2\pi.$  This information is referred to as 
the {\it persistence for the circle valued map $f$}.

When $f$ is tame, this is again completely determined up to 
coherent isomorphisms by a finite collection of invariants.
However, unlike sublevel and level persistence for real valued
maps, the invariants
include structures other than bar codes called {\em Jordan cells}.
Specifically, for any $r= 0,1, \cdots, \dim(X)$ 
we have two types of invariants:
\begin{itemize}
\item {\em bar codes}: intervals 
with ends $s,s'$
$0<s\leq 2\pi, \ s\leq s' < \infty$,
that are closed or open at $s$ or $s'$, 
precisely of one of the forms $[s,s'], (s, s'], [s,s')$, and $(s,s').$
These intervals can be geometerized as  ``spirals"
with equations in (\ref{E1}). For any interval 
$\{s,s'\}$ the spiral is the plane curve 
(see Figure~\ref{spiral} in section \ref{S4})

\begin{equation}\label{E1}
\begin{array}{ccc}
\begin{array}{c}
x(\theta)= \left(\theta + 1-s\right) \cos \theta  \\
y(\theta)=\left(\theta+1-s\right) \sin \theta 
\end{array} & with  &
 \theta\in \{s, s'\}.
\end{array}
\end {equation}


\item {\em Jordan cells}.
A Jordan cell is a pair $(\lambda, k),$  
$\lambda \in \overline{\kappa} \setminus 0, \ \  k \in \mathbb Z_{>0},$ where $\overline{\kappa}$ 
denotes the algebraic closure of the field $\kappa.$
It corresponds to a $k\times k$ matrix of the form
\begin{equation}
\begin{pmatrix}
\lambda&1        & 0\dots    &0\\
0 & \lambda & 1\dots & 0\\
\vdots\\
0             & \dots&\lambda &1\\
0             & \dots&0              &\lambda 
\end{pmatrix}.
\label{eqnmat}
\end{equation}

\item{$r$-invariants.}  Given a tame map $f: X\rightarrow \mathbb{S}^1$,
the collection of bar codes and Jordan cells for 
each dimension $r\in \{ 0,1,2,\cdots \dim X\}$ 
constitute the {\it $r$-invariants} 
of the map $f.$  
\end{itemize}

We will define all of the above items 
in the next section using quiver representations.

The bar codes for $f$ can be inferred 
from $\tilde f: \tilde X_{[a,b]}\to \mathbb R$ with 
$[a,b]$ being any large enough interval.
Specifically, the bar codes of $f:X\to {\mathbb S}^1$ 
are among the ones of $\tilde f: \tilde X_{[a,b]}\to \mathbb R$ 
for $(b-a)$ being at most
$\sup_\theta \dim H_r(X_\theta).$ 

The Jordan cells can not be derived  
from $\tilde f:\tilde X\to \mathbb R$ or any 
of its truncations $\tilde f:\tilde X_{[a,b]} \to \mathbb R$ unless additional information, like the deck transformation of $\tilde X,$ is provided.
The end points of any bar code for $f$ 
correspond to critical angles, that is,
$s$ and $s'\pmod{2\pi}$
of a bar code interval $\{s,s'\}$ 
are critical angles for $f$. One can recover the following information
from the bar codes and Jordan cells:
\begin{enumerate}
\item The Betti numbers of each fiber,
\item The Betti numbers of the source space $X$, and
\item  The dimension of the kernel and the image of the linear map induced in homology by the inclusion $X_\theta \subset X$
 as well as other additional topological invariants not discussed here~\cite {BH12}.
\end{enumerate}


Theorems~\ref{T1} and \ref{T2} make the above statement precise.
Let $B$ be a bar code 
described by a spiral in (\ref {E1})
and $\theta$ be any angle. Let 
$n_\theta(B)$ denote the cardinality of the intersection of 
the spiral with
the ray originating at the origin and making 
an angle $\theta$ with the $x$-axis. For the 
Jordan cell $J=(\lambda, k)$, let  $n(J)=k$ and $\lambda(J)= \lambda.$
Furthermore, let $\BB_r$ and ${\JJ}_r$ denote the set of bar codes 
and Jordan cells for $r$-dimensional homology. We have the following results.

\begin{theorem}
$ \dim H_r (X_\theta)=\sum _{B\in \mathcal B_r}  n_\theta (B) + \sum_{J\in \mathcal J_r} n(J).
$
\label {T1}
\end{theorem}   

\begin{theorem} 
$ \dim H_r (X)=\#\{B\in \BB_r| {\rm both\  ends\  closed}\} 
+ \#\{B\in \BB_{r-1} | {\rm both\  ends\  open}\} + 
 \# \{J\in \JJ_{r}| \lambda(J)=1\} +\# \{J\in \JJ_{r-1}| \lambda(J)=1\}.
$ 
\label {T2}
\end{theorem}

Using the same arguments as in the proof of the
above Theorems one can derive:
\begin{proposition}\label{P33}
$\dim \rm{img}(H_r(X_\theta)\to H_r(X))= \#\{B\in \BB_r| n_\theta(B)\ne 0 \ {\rm and\ both\  ends\  closed}\}  + \# \{J\in \JJ_{r}| \lambda=1\}$
\end{proposition}

A real valued tame map $f: X\rightarrow \mathbb{R}$ can be
regarded as a circle valued tame map $f':X\rightarrow \mathbb{S}^1$
by identifying $\mathbb{R}$ to $(0,2\pi)$ with critical values
$t_1,\cdots, t_m$ becoming the critical angles 
$\theta_1,\cdots,\theta_m$ where $\theta_i=2\arctan t_i + \pi$. 
The map $f'$ in this case will not have any 
Jordan cells and the bar codes will be the same as level persistence bar codes. We have the following corollary:

\begin{corollary} 
$ \dim H_r (X_\theta)=\sum _{B\in \mathcal B_r}  n_\theta (B) $ and

$\dim H_r (X)=\#\{B\in \BB_r| \rm{ both\  ends\  closed}\} 
+ \#\{B\in \BB_{r-1} | \rm{ both\  ends\  open}\}.
$
\label {C1}
\end{corollary}

Theorem~\ref{T1} is quite intuitive and is in analogy
with the derived results for sublevel and level persistence~\cite{CSD09,ZC05}.
Theorem~\ref{T2} is more subtle. Its counterpart
for real valued function (Corollary~\ref{C1}) has not yet
appeared in the literature though a related result for homology
of source space can be derived from extended persistence~\cite{CEH09}. 
The proofs of these results require the definition of the bar codes 
and Jordan cells which appear in the next section. 
The proofs are sketched in section \ref {S7}. 

The Questions Q1-Q3 can be answered 
using the bar codes. The question Q4 about returned homology
can be answered using the bar codes and Jordan cells.
\vskip .1in


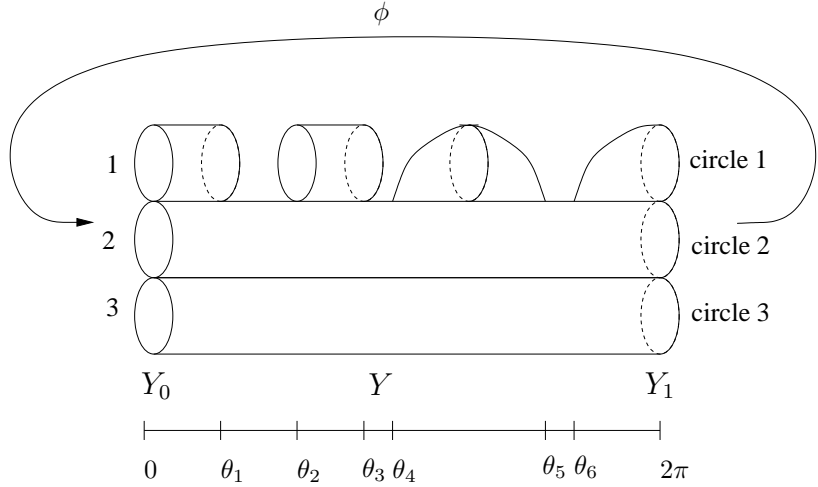
\begin{figure}[h]
\centerline{\input{cv-ex.pstex_t}}
~\\
\begin{tabular}{c|c}
map $\phi$ & $r$-invariants \\ \hline
\begin{tabular}{l}
circle 1: 1 times around 1,-3 times around 2, -2 times around 3\\
circle 2: 1 times around 1, 4 times around 2, 1 time around 3\\
circle 3: 2 times around 1, 2 times around 2, 2 times around 3  \\
\end{tabular} & 
\begin{tabular}{c|c|c}
dimension & bar codes & Jordan cells \\ \hline
0  &  & $(1,1)$ \\ \hline
   & $(\theta_6, \theta_1+2\pi]$ & $(3,2)$ \\
1  & $[\theta_2,\theta_3]$ & \\
   & $(\theta_4,\theta_5)$ &  \\
\end{tabular} \\ \hline 
\end{tabular}
\caption{Example of $r$-invariants for a circle valued map}
\label{cv-ex}
\end{figure}

Figure~\ref{cv-ex} indicates a tame map $f:X\to {\mathbb S}^1$ and the 
corresponding invariants, bar codes, and Jordan cells.
The space $X$ is obtained from $Y$ in the figure by identifying its  
right end $Y_1$ (a union of three circles) to the left end 
$Y_0$ (again a union of three circles ) following 
the map $\phi:Y_1\to Y_0.$  The map $f:X\to {\mathbb S}^1$ is induced by 
the projection of $Y$ on the interval $[0,2\pi].$
We have $H_1(Y_1)=H_1(Y_0)=\kappa\oplus\kappa\oplus\kappa$
and $\phi$ induces a linear  map in $1$-homology represented by the matrix  
\footnote {Each circle is oriented counterclockwise and represents a 
$1$-dimensional homology class; ``$k$ times ($-k$ times) 
around the circle"  means " going around  $k$ times  
counter clockwise (clockwise respectively)".} 
\begin{equation*}
\begin{pmatrix}
1&1&2\\
-3&4&2\\
-2&1&2            
\end{pmatrix}.
\end{equation*}
The first generator (circle 1) of $H_1(\tilde X_{2\pi})$ is dead in $H_1(\tilde X_{[\theta,2\pi]})$ for $\theta\leq \theta_6$ but not for $\theta\in (\theta_6, 2\pi]$ and is  detected in $H_1(\tilde X_{2\pi+\theta})$ for $0\leq \theta\leq \theta_1$ but not for $\theta>\theta_1.$ It generates a bar code $(\theta _6,  2 \pi +\theta_1].$ The other two (circle 2 and 3) 
never die and provide a Jordan cell $(3,2).$ In Appendix 
we show how our algorithm can be used to compute the bar codes and 
Jordan cells for the above example.


\section {Representation theory and its invariants}\label{S4}
The invariants for the circle valued map are
derived from the representation theory of quivers. 
The quivers are directed graphs. The representation
theory of simple quivers such
as paths with directed edges was described by Gabriel~\cite{G72} 
and is at the heart 
of the derivation of the invariants for 
zigzag and then level persistence in~\cite{CSD09}. For circle valued
maps, one needs  representation theory for
circle graphs with directed edges.
This theory appears in the work of Nazarova~\cite {N73}, and
Donovan and  Ruth-Freislich~\cite {DF73}. 

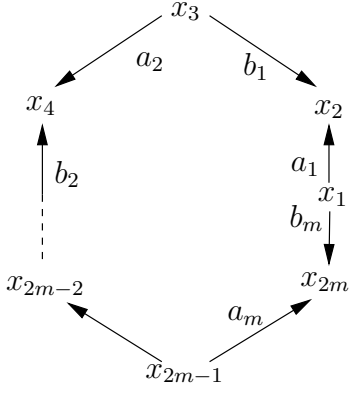
\parpic[l]{
\begin{picture}(0,0)%
\includegraphics{quiver.pstex}%
\end{picture}%
\setlength{\unitlength}{3947sp}%
\begingroup\makeatletter\ifx\SetFigFont\undefined%
\gdef\SetFigFont#1#2#3#4#5{%
  \reset@font\fontsize{#1}{#2pt}%
  \fontfamily{#3}\fontseries{#4}\fontshape{#5}%
  \selectfont}%
\fi\endgroup%
\begin{picture}(2622,2489)(8551,-2769)
\put(10801,-1051){\makebox(0,0)[b]{\smash{{\SetFigFont{12}{14.4}{\rmdefault}{\mddefault}{\updefault}{\color[rgb]{0,0,0}$x_2$}%
}}}}
\put(10261,-811){\makebox(0,0)[lb]{\smash{{\SetFigFont{12}{14.4}{\rmdefault}{\mddefault}{\updefault}{\color[rgb]{0,0,0}$b_1$}%
}}}}
\put(9586,-751){\makebox(0,0)[lb]{\smash{{\SetFigFont{12}{14.4}{\rmdefault}{\mddefault}{\updefault}{\color[rgb]{0,0,0}$a_2$}%
}}}}
\put(9076,-1486){\makebox(0,0)[lb]{\smash{{\SetFigFont{12}{14.4}{\rmdefault}{\mddefault}{\updefault}{\color[rgb]{0,0,0}$b_2$}%
}}}}
\put(9901,-436){\makebox(0,0)[b]{\smash{{\SetFigFont{12}{14.4}{\rmdefault}{\mddefault}{\updefault}{\color[rgb]{0,0,0}$x_3$}%
}}}}
\put(9896,-2711){\makebox(0,0)[b]{\smash{{\SetFigFont{12}{14.4}{\rmdefault}{\mddefault}{\updefault}{\color[rgb]{0,0,0}$x_{2m-1}$}%
}}}}
\put(9021,-2161){\makebox(0,0)[b]{\smash{{\SetFigFont{12}{14.4}{\rmdefault}{\mddefault}{\updefault}{\color[rgb]{0,0,0}$x_{2m-2}$}%
}}}}
\put(8991,-1031){\makebox(0,0)[b]{\smash{{\SetFigFont{12}{14.4}{\rmdefault}{\mddefault}{\updefault}{\color[rgb]{0,0,0}$x_4$}%
}}}}
\put(10651,-1411){\makebox(0,0)[b]{\smash{{\SetFigFont{12}{14.4}{\rmdefault}{\mddefault}{\updefault}{\color[rgb]{0,0,0}$a_1$}%
}}}}
\put(10661,-1761){\makebox(0,0)[b]{\smash{{\SetFigFont{12}{14.4}{\rmdefault}{\mddefault}{\updefault}{\color[rgb]{0,0,0}$b_m$}%
}}}}
\put(10281,-2351){\makebox(0,0)[b]{\smash{{\SetFigFont{12}{14.4}{\rmdefault}{\mddefault}{\updefault}{\color[rgb]{0,0,0}$a_m$}%
}}}}
\put(10781,-2136){\makebox(0,0)[b]{\smash{{\SetFigFont{12}{14.4}{\rmdefault}{\mddefault}{\updefault}{\color[rgb]{0,0,0}$x_{2m}$}%
}}}}
\put(10826,-1601){\makebox(0,0)[b]{\smash{{\SetFigFont{12}{14.4}{\rmdefault}{\mddefault}{\updefault}{\color[rgb]{0,0,0}$x_1$}%
}}}}
\end{picture}%

}
Let $G_{2m}$ be a directed graph with $2m$ vertices, 
$x_1, x_1, \cdots x_{2m}$. Its underlying
undirected graph is a simple cycle. The directed edges in $G_{2m}$
are of two types: {\em forward}
$a_i: x_{2i-1}\to x_{2i}$, $1\leq i\leq m$, and 
{\em backward} $b_i: x_{2i+1}\to x_{2i}$, $1\leq i\leq m-1$, 
$b_m:x_1\to x_{2m}$.

We think of this graph as being located
on the unit circle centered at the origin $o$ in the plane.

A representation $\RR$ on $G_{2m}$ is an
assignment of a vector space $V_x$ to each vertex $x$
and a linear map $\ell_e: V_{x}\rightarrow V_{y}$ for each oriented
edge $e=\{x,y\}$. Two representations $\RR$ and $\RR'$
are isomorphic if for each vertex $x$ there 
exists an isomorphism from the vector space $V_x$ of 
$\RR$ to the vector space $V'_x$ of  $\RR'$,
and these isomorphisms commute with
the linear maps $V_{x}\rightarrow V_{y}$ and
$V'_{x}\rightarrow V'_{y}$.  
A {\em non-trivial} representation assigns at least one
vector space which is not zero-dimensional.
A representation is {\em indecomposable}
if it is not isomorphic to the sum of two nontrivial representations.

Given two representations $\rho$ and $\rho'$, their sum $\rho\oplus \rho'$
is a representation whose vector spaces are the direct sums 
$V_x\oplus V'_x$ related by
linear maps that are the direct sums $\ell_e\oplus \ell'_e.$
It is not hard to observe that each representation 
has a decomposition as a sum of indecomposable representations 
unique up to isomorphisms.

We provide
a description of indecomposable representations of the quiver $G_{2m}.$
For any triple of integers $\{i,j,k\}$, $1 \leq i,j \leq m,$ $k\geq 0$,
one may have any of the 
four representations, 
$\RR^I([i, j];k), \ $  $\RR^I((i,j];k), $ \newline $\RR^I([i,j);k) \ $,
and $\RR^I((i,j);k) \ $ defined below.  For any Jordan cell $(\lambda, k)$ one has the representation $\rho^J(\lambda, k)$ defined below. 
The exponents $I$ and $J$ indicate that these representations 
are associated with a bar code (interval) or a Jordan cell
respectively and hence we call them bar code and Jordan cell
representations.

\begin{itemize}
\item Bar code representation $\RR^I(\{i, j\};k)$:   
Suppose that the  evenly indexed vertices $\{x_2, x_4, \cdots x_{2m}\}$ 
of $G_{2m}$ which are the targets of the directed arrows  
correspond to the angles 
$0< s_1< s_2 <\cdots <s_{m} \leq2\pi.$ 
Draw the spiral curve given by (\ref{E1})
for the interval $\{s_{i}, s_{j}+2k\pi\}$; refer to 
Figure~\ref{spiral}.

For each $x_i$, let $\{e_i^1,e_i^2,\cdots\}$ denote the 
ordered set (possibly empty) of intersection points of the
ray $ox_i$ with the spiral. 
While considering these intersections, it is important to realize that 
the point $(x(s_i), y(s_i))$ (resp. $(x(s_j+2k\pi), y(s_j+2k\pi))$) 
does not belong to the spiral (\ref{E1})
if  $\{i,j\}$  
is open at $i$ (resp. $j$). 
For example, in Figure~\ref{spiral},
the last circle on the ray $ox_{2j}$ 
is not on the spiral since $[i,j)$ in $\rho^I([i,j);2)$ is open at right. 

Let $V_{x_i}$ denote the vector space
generated by the base $\{e_i^1,e_i^2,\cdots\}$. 
Furthermore, let
$\alpha_i: V_{x_{2i-1}}\rightarrow V_{x_{2i}}$ and 
$\beta_i: V_{x_{2i+1}}\rightarrow V_{x_{2i}}$ be the linear maps
defined on bases and extended by linearity as follows: assign
the vector $e^h_{2i}\in V_{x_i}$ to $e^{\ell}_{2i\pm 1}$
if $e^h_{2i}$ is an 
adjacent intersection point 
to the points $e^{\ell}_{2i\pm 1}$ on the spiral. If $e^h_{2i}$ does
not exist, assign zero to $e^{\ell}_{2i\pm 1}$.
If $e^{\ell}_{2i\pm 1}$ do not go to zero,
$h$ has to be $l$, $l-1$, or $l+1$.
The construction above provides a representation on 
$G_{2m}$ which is indecomposable. 
Once the angles $s_i$ are associated to the vertices $x_{2i}$ one 
can also think of these representations $\rho^I(\{i,j\};k)$ as the bar codes  
$[s_{i}, s_{j}+2k\pi]$, $(s_{i}, s_{j}+2k\pi],$ $[s_{i}, s_{j}+2k\pi)$,
and $(s_{i}, s_{j}+2k\pi)$.
\begin{figure}[h]
\centerline{\input{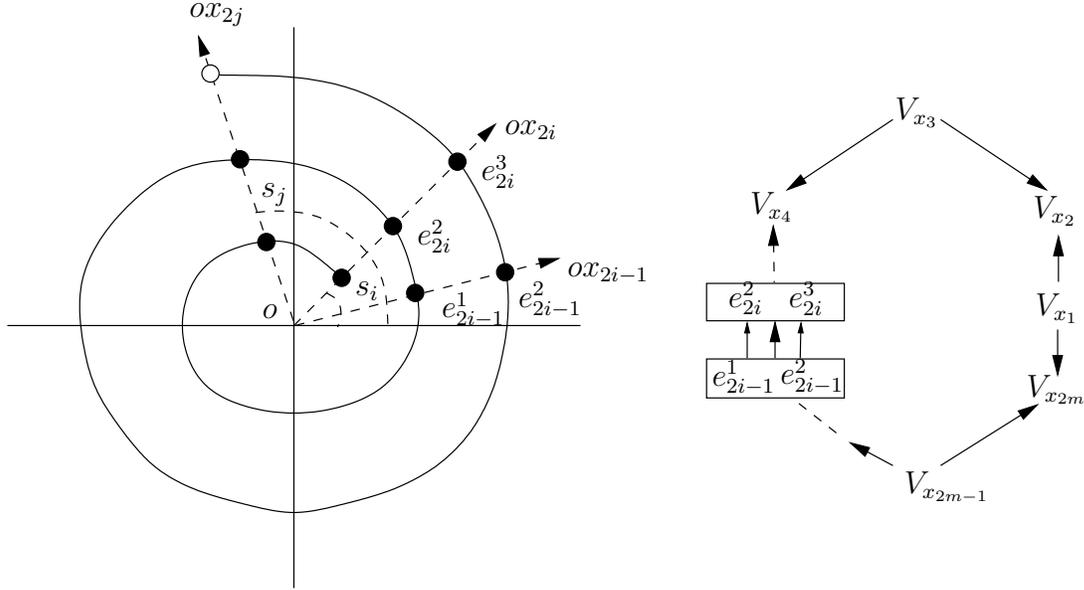}}
\caption{The spiral for $[s_i,s_j+4\pi)$.}
\label{spiral}
\end{figure} 


\item Jordan cell representation $\RR^J(\lambda,k)$:
 Assign the vector space 
with the base $\{e_1,e_2,\cdots,e_k\}$ to each $x_i$ and take all linear
maps $\alpha_i$ but one (say $\alpha_1$) and $\beta_i$  the
identity. The linear map $\alpha_1$ is given by the Jordan matrix defined by 
$(\lambda,k)$ in (2).  
Again this representation is indecomposable.
\end{itemize}

It follows from the work of~\cite{DF73,N73}
that when $\kappa$ is algebraically closed\footnote {when $\kappa$ is not algebraically closed Jordan cells  have to be replaced by  conjugate classes of  indecomposable (not conjugated to a direct sum of matrices) matrices with entries in $\kappa.$ }, 
the bar code and Jordan cell representations are all and only
indecomposable representations of the quiver $G_{2m}.$   
The collection of all bar code and Jordan cell representations of
a representation $\rho$ constitutes its {\it invariants}.

Now, consider the representation  $\rho$ on the graph  $G_{2m}$ given 
by the vector spaces $V_{2i-1}:=V_{x_{2i- 1}}, V_{2i}:=V_{x_{2i}}$ 
and the linear maps $\alpha_i$ and $\beta_i.$ 
To such a representation $\rho$, we associate a map 
$M_\rho:  \bigoplus_{1\leq i\leq m} V_{2i- 1} \to   
\bigoplus _{1\leq i\leq m}V_{2i}$ 
which is represented by a block matrix also denoted as $M_{\rho}$: 

\begin{equation*}
\begin{pmatrix}

\alpha_1& -\beta_1 &      0       &\dots  & \dots &0\\
         0    & \alpha_2 &-\beta_2&\dots   &  \dots &0\\
     \vdots&\vdots       &\vdots     &\vdots &\vdots & \\
      0        &\dots        &\dots       &\dots    &\dots \alpha_{m-1}&-\beta_{m-1}\\
      -\beta_m&\dots    &\dots       &\dots   & \dots &\alpha_m
         \end{pmatrix} 
         \end{equation*}

For this representation we define its dimension characteristic
as the $2m$-tuple 
of positive integers 
$$\dim (\rho)= (n_1,r_1 \cdots  n_m, r_m)$$ 
with $n_i= \dim V_{x_{2i-1}}$ and $r_i=\dim V_{x_{2i}}$ and  
denote by $\ker (\rho):=  \ker M_\rho$ and 
$\coker (\rho)= \coker M_\rho.$ For the sum of 
two such representations $\rho= \rho_1 \oplus \rho_2$ we have:

\begin{proposition} \label {PR}
~
\begin{enumerate}
\item $\dim(\rho_1\oplus \rho_2)= \dim(\rho_1) + \dim(\rho_2)$,  
\item $\dim \ker (\rho_1\oplus \rho_2)=\dim \ker (\rho_1) + \dim \ker (\rho_2)$,
\item $\dim \coker(\rho_1\oplus \rho_2)= \dim \coker (\rho_1) + 
\dim \coker (\rho_2)$.
\end{enumerate}
\end{proposition}
The description of a bar code representation permits explicit calculations.

\begin{proposition}\label {AP1}
~
\begin{enumerate}
\item If $i\leq j$ then 
\begin{enumerate}
\item $\dim \rho^I([i,j]; k)$ is given by: 

$n_l=k+1$ if  $(i+1)\leq l\leq j$ and $k$ otherwise, 

$ r_l= k+1$ if $i\leq l\leq j$ and $k$ otherwise 
\item $\dim \rho^I((i,j]; k)$ is given by: 

$n_l=k+1$ if  $(i+1)\leq l\leq j$ and  $k$ otherwise, 

$r_l= k+1$ if $(i +1)\leq l\leq j$ and $k$ otherwise, 
\item $\dim \rho^I([i,j); k)$ is given by: 

$n_l=k+1$ if  
$(i+1)\leq l\leq j$ and $k$ otherwise,  

$ r_l= k+1$   if  $i\leq l\leq (j-1)$ and $k$ otherwise,   
\item $\dim \rho^I((i,j); k)$ is given by:

 $n_l=k+1$ if  
$(i+1)\leq l\leq j$ and $k$ otherwise,  

$r_l= k+1$  if  $(i +1)\leq l\leq( j-1)$ and $k$ otherwise

\end{enumerate}
\item If $i> j$ then similar statements hold.
\begin{enumerate}
\item $\dim \rho^I([i,j]; k)$ is given by:

 $n_l=k$ if  $(j+1)\leq l\leq i$ and $k+1$ otherwise; 
 
 $r_l=k$ if   $(j+1)\leq l\leq (i-1)j$ and  $k+1$ otherwise
\item $\dim \rho^I((i,j]; k)$ is given by:

 $n_l=k$ if  $(j+1)\leq l\leq i$ and  
 $k+1 $ otherwise.  
 
 $r_l=k$ if  
 $(j+1)\leq l\leq i$ and $k +1$ otherwise,
\item $\dim \rho^I([i,j); k)$ is given by:

 $n_l=k$ if  
$(j+1)\leq l\leq i$ and  $k+1$ otherwise;

 $r_l=k$ if  
$ j\leq l \leq (i-1)$  and $k+1 $ otherwise,
\item $\dim \rho^I((i,j); k)$ is given by:

 $n_l=k$  if $(j+1)\leq l\leq i$  and $k+1$ otherwise; 
 
 $r_l=k$ if  
$j\leq l\leq i$  and $k+1$ otherwise.  
\end{enumerate}
\item $\dim\rho^J(\lambda,k)$ is given by $n_i=r_i=k$
\end{enumerate}
\end{proposition}


\begin{proposition}\label{AP2}
~
\begin{enumerate} 
\item $\dim \ker \rho^I([i,j]; k)=0$, $\dim\coker \rho^I([i,j]; k)=1$,
\item $\dim \ker \rho^I([i,j); k)=0$, $\dim \coker \rho ^I([i,j); k)=0$,
\item $\dim \ker \rho^I ((i,j]; k)=0$, $\dim \coker \rho^I((i,j]; k)=0$,
\item $\dim \ker \rho^I((i,j); k)=1$, $\dim \coker \rho^I((i,j); k)=0$,
\item $\dim \ker \rho^{J}(\lambda,k) =0$ (\rm{resp.} 1) 
if $\lambda {\ne 1}$ (\rm{resp.} 1),
\item $\dim \coker \rho^{J}(\lambda,k)=0$ (\rm{resp.} 1) 
if $\lambda {\ne 1}$ (\rm{resp.} 1).
\end{enumerate}
\end{proposition}

\begin{obs}  \label{O1}
A representation  $\RR$ has no  indecomposable 
components of type $\RR^I$ in its decomposition
iff all linear maps $\alpha_i'$s and $\beta_i'$s are isomorphisms. 
For such a representation, starting with an index $i,$ 
consider the linear isomorphism 
$$T_i=\beta^{-1}_{i}\cdot \alpha_{i}\cdot \beta^{-1}_{i-1}\cdot \alpha_{i-1}\cdots \beta^{-1}_{2}\cdot \alpha_{2}\cdot \beta^{-1}_{1}\cdot \alpha_{1}\cdot \beta^{-1}_m\cdot \alpha_m\cdot \beta^{-1}_{m-1}\cdot \alpha_{m-1} \cdots \beta^{-1}_{i+1}\cdot \alpha_{i+1}.$$ 
The Jordan canonical form~\cite{DS58} 
of the isomorphism $T_i$ is independent of $i$ 
and is a block diagonal matrix with the diagonal consisting of 
Jordan cells $(\lambda, k)$s.
Clearly, $\RR$ is the direct sum of 
$\RR^J(\lambda,k)$s, the Jordan cell representations of $\rho$.
\end{obs}


\begin{definition}[$r$-invariants.]
Let $f$ be a circle valued tame map defined on a topological space $X$.
For $f$ with $m$ critical angles $0< s_1< s_2,\cdots s_m\leq 2\pi$,
consider the quiver 
$G_{2m}$ with  the vertices  $x_{2i}$ identified with 
the angles $ s_i$ and the vertices $x_{2i-1}$ identified with the 
angles $t_i$ that satisfy $0<t_1<s_1<t_2 < s_2, \cdots t_m<s_m.$ 

For any $r$, consider the representation  $\rho_r$ of $G_{2m}$
with $V_{x_i}=H_r(X_{x_i})$ and the linear maps $\alpha_i$s and $\beta_i$s
induced in the $r$-homology by maps 
$a_i: X_{x_{2i-1}} \to X_{x_{2i}}$ 
and $b_i: X_{x_{2i+1}} \to X_{x_{2i}}$ described in section \ref{S1}. 
The bar code and Jordan cell representations 
of $\rho_r$ are independent of the choice of $t_i$s  
and are collectively referred as the  {\it $r$-invariants} of $f.$
\end{definition}

\section{ Proof of the main results}\label{S7}

The Figure~\ref{cv-ex}
and the bar codes listed below suggest why a semi-closed 
(one end open and the other closed) bar code does not contribute 
to the homology of the total space $X$ and
why a closed bar code (both ends closed) in ${\cal B}_r$ 
contributes one unit while an open (both end open) bar 
code  in ${\cal B}_{r-1}$ contributes one unit to the 
$H_r(X)$. 
Indeed, in our example, a  semi-closed bar code in
${\cal B}_1$ adds to the total 
space a cone over ${\mathbb S}^1,$ which is a contractible space. 
It gets glued to the total
space along a generator of the cone (a segment connecting the apex to
${\mathbb S}^1$), again a contractible space. 
A closed bar code in ${\cal B}_1$  adds a cylinder of 
${\mathbb S}^1$ whose $H_1$ has dimension $1$.
It gets glued to the total space along a generator of the cylinder (a segment
connecting the same point on the two copies
of ${\mathbb S}^1$), again a contractible space. An open bar code
in ${\cal B}_1$ adds the suspension over ${\mathbb S}^1$, 
topologically a $2$-sphere which gets glued along a meridian, 
a contractible space. This contributes a dimension to $H_2$. 

The lack of contribution of a Jordan cell with $\lambda \ne 1$ 
as well as  the contribution of one unit of a 
Jordan cell in ${\cal J}_r$ with $\lambda=1$ to both $r$ and $r+1$ 
dimensional homology of the total space
should not be a surprise for the reader familiar with the 
calculation of the homology of mapping torus. 

Below we explain rigorously but 
schematically the arguments  for the proof of Theorems~\ref{T1}, 
~\ref{T2}, and Corollary~\ref{C1}.

The proof of Theorem~\ref{T1} is a consequence 
of Propositions~\ref{PR} and ~\ref{AP1}.
The proof of Theorem~\ref{T2} proceeds along the following lines.

First observe that, up to homotopy, the space $X$  can be 
regarded as the iterated mapping torus $\mathcal T$ described below. 
Consider the collection of spaces and continuous maps:
\begin{eqnarray*}
X_m =X_0 \stackrel{b_0=b_m}{\longleftarrow} R_1
\stackrel{a_1}{\longrightarrow}X_1
\stackrel{b_1}{\longleftarrow}R_2
\stackrel{a_2}{\longrightarrow}X_2 \cdots X_{m-1}
\stackrel{b_{m-1}}{\longleftarrow}R_m 
\stackrel{a_m}{\longrightarrow}X_m
\end{eqnarray*}
 with  $R_i:=X_{t_i}$ and $X_i:= X_{s_i}$
and denote by $\mathcal T=T(\alpha_1 \cdots \alpha_m; 
\beta_1\cdots\beta_m)$ the space obtained from the disjoint union 
$$(\bigsqcup_{1\leq i\leq m} R_i\times [0,1])
\sqcup(\bigsqcup_{1\leq i\leq m} X_i)
$$
by identifying $R_i\times \{1\}$ to $X_i$ by $\alpha_i$  
and $R_i\times \{0\}$ to $X_{i-1}$ by $\beta_{i-1}.$ 
Denote by $f^{\mathcal T} :\mathcal T\to [0,m]$ where 
$f^{\mathcal T}: R_i\times [0,1]\to [i-1,i]$ is the projection on
$[0,1]$ followed by the translation of $[0,1]$ to 
$[i-1,i]$. This  map is a homotopical reconstruction of 
$f: X \to {\mathbb S}^1$ provided that,  with the choice of angles  $t_i, s_i$ and maps $a_i$ $b_i$ described in section \ref {S1}, $X_i:= f^{-1}(s_i) , R_i:= f^{-1}(t_i)$.

Let $\mathcal P'$ denote the space obtained from the disjoint union 
$$(\bigsqcup_{1\leq i\leq m} R_i\times (\epsilon,1])\sqcup(\bigsqcup_{1\leq i\leq m} X_i)$$
by identifying $R_i\times \{1\}$ to $X_i$ by $\alpha_i,$ and   
$\mathcal P''$ denote the space obtained from the disjoint union 
$$(\bigsqcup_{1\leq i\leq m} R_i\times [0,1-\epsilon)\sqcup(\bigsqcup_{1\leq i\leq m} X_i)$$
by identifying $R_i\times \{0\}$ to $X_{i-1}$ by $\beta_{i-1}.$

Let $\mathcal R= \bigsqcup_{1\leq i\leq m} R_i$ and $\mathcal X= \bigsqcup_{1\leq i\leq m} X_i.$ Then, one has:
\begin{enumerate}
\item $\mathcal T= \mathcal P' \cup \mathcal P''$,
\item $\mathcal P' \cap \mathcal P''= 
(\bigsqcup_{1\leq i\leq m} R_i\times (\epsilon,1-\epsilon))\sqcup \mathcal X$,
and
\item the inclusions $(\bigsqcup_{1\leq i\leq m}  
R_i\times \{1/2\})\sqcup \mathcal X\subset \mathcal P'\cap \mathcal P''$
as well as the obvious inclusions 
$\mathcal X\subset \mathcal P'\rm{and} \  \mathcal X\subset  \mathcal P''$
are homotopy equivalences.
\end{enumerate}
The Mayer-Vietoris long exact sequence leads to the diagram 

\xymatrix{
&&                                                       & H_r(\mathcal R)\ar[r]^{M_{\rho_r}}                                                  & H_r(\mathcal X)\ar[rd]                                               &&\\
&\cdots\ar[r]&H_{r+1}(\mathcal T)\ar[ur]
\ar[r]^{\partial_{r+1}}&H_r(\mathcal R)\oplus H_r(\mathcal X)\ar[u]^{pr_1}\ar[r]^N  &H_r(\mathcal X)\oplus H_r(\mathcal X)\ar[u]^{(Id, -Id)}{\ar[r]^{(i^r, -i^r)}} &H_r(\mathcal T)\ar[r]&\\
&&                                                       &H_r(\mathcal X)\ar[u]^{in_2}\ar[r]^{Id}                                           &H_r(\mathcal X)\ar[u]^{\Delta}                                                 &&}

Here $\Delta$ denotes the diagonal, $in_2$ the inclusion on the second component, $pr_1$ the projection on the first component, 
$i^r$ the linear map induced in homology by 
the inclusion $\mathcal X\subset \mathcal T$, and 
$M_{\rho_r}$ the map given by the matrix 

\begin{equation}
\begin{pmatrix}

\alpha^r_1& -\beta^r_1 &      0       &\dots  & \dots &0\\
         0    & \alpha^r_2 &-\beta^r_2&\dots   &  \dots &0\\
     \vdots&\vdots       &\vdots     &\vdots &\vdots & \\
      0        &\dots        &\dots       &\dots    &\dots \alpha^r_{m-1}&-\beta^r_{m-1}\\
      -\beta^r_m&\dots    &\dots       &\dots   & \dots &\alpha^r_m
         \end{pmatrix} 
\label{main-matrix}
\end{equation}
with $\alpha^r_i : H_r(R_i)\to H_r(X_i)$ and $\beta^r_i: H_r(R_{i+1})\to H_r(X_i)$ induced by the maps $\alpha_i$ and $\beta_i,$
and $N$ defined by 
\begin{equation*}
\begin{pmatrix}
&\alpha^r & Id\\
&-\beta^r &Id
\end{pmatrix}
\end{equation*}
where  $\alpha^r$ and $\beta^r$ are the matrices 
\begin{equation*}
\begin{pmatrix}
\alpha^r_1&  0       &\dots  & \dots &0\\
         0    & \alpha^r_2 &\dots   &  \dots &0\\
     \vdots&\vdots        &\vdots &\vdots & \vdots \\
      0           &0       &\dots    &0 &\alpha^r_{m-1} 
\end{pmatrix}
\end{equation*}
\vskip .1in

\begin{equation*}
\begin{pmatrix}
0&  \beta^r_1 &      0        & \dots &0\\
         0    &0 & \beta^r_2& \dots           &0\\
     \vdots&\vdots       &\vdots &\vdots & \vdots \\
      0        &\dots        &\dots        & 0    &\beta^r_{m-1}\\
      \beta^r_m&0         &\dots   & 0&0
\end{pmatrix}.
\end{equation*}

From the diagram above we retain only the long exact sequence
\begin{equation}\label{MV}
\cdots\to H_r(\mathcal R)\xrightarrow{M_{\rho_r}}H_r(\mathcal X)\to H_r(\mathcal T)\to H_{r-1}(\mathcal R)\xrightarrow{M_{\rho_{r-1}}}H_{r-1}(\mathcal X)\to\cdots
\end{equation}
 from which we derive 
 the short exact sequence 
 
 \begin{equation}\label{EE13}
0\to\coker (\rho_r)\to H_r(\mathcal T)\to\ker (\rho_{r-1})\to0
\end{equation}

and then 
\begin{equation}
H_r(\mathcal T)= \rm{co ker} \rho_r \oplus \ker \rho_{r-1}
\label{eqf}
\end{equation}

Theorem \ref{T2}  follows from Propositions \ref{PR}, \ref{AP2} and
the equation~(\ref{eqf}) above.
A specified decomposition of $\rho_r$ and $\rho_{r-1}$ 
into indecomposable representations  and a splitting 
in the sequence~ (\ref{EE13}) provide specified elements in $H_r(X_\theta)$ and $H_r(\mathcal T)$ which can be compared. This leads to  the verification of Proposition  \ref{P33}.

\section {Algorithm}\label{S5}
Given a circle valued tame map $f:X\rightarrow \mathbb{S}^1$, we now
present an algorithm to compute the bar codes and the Jordan cells when $X$ 
is a finite simplicial complex, and $f$ is generic and linear. This makes the map tame.
Genericity means that $f$ is injective on vertices.  
To explain linearity we recall that, for 
any simplex $\sigma \in X$, the restriction  
$f|_{\sigma}$ admits liftings $\hat
f:\sigma \to \mathbb R$, i.e. $\hat f$ is a 
continuous map which satisfies $p\cdot
\hat f= f|_\sigma.$
The map $f: X\to \mathbb S^1$ is called  {\em linear} 
if for any simplex $\sigma $, at least
one  of the liftings (and then any other) is  linear.

Our algorithm takes the simplicial complex $X$ equipped with the map $f$
as input and, for any $r$, computes the matrix $M_{\rho_r}$ 
of the representation
$\rho_r$ for $f$. This requires recognizing the critical values
$s_1, s_2, \cdots s_{m}\in \mathbb S^1$ of $f$, and 
for conveniently 
chosen regular values $t_1, t_2, \cdots t_{m}\in \mathbb S^1$, 
determining the vector spaces $V_{2i-1}=H_r(X_{t_i}), 
V_{2i}=H_r(X_{s_i})$ with
the linear maps $\alpha_i$ and $\beta_i$ as matrices.
We consider the block matrix 
$M_{\rho_r}:  \bigoplus_{1\leq i\leq m} V_{2i-1} \to
\bigoplus _{1\leq i\leq m}V_{2i}$ described in the previous section.

We compute the bar codes from the  block matrix $M_{\rho_r}$ first, and
then the Jordan cells.
The algorithm 
consists of three steps. We describe the first and second steps in sufficient
details. The third step is a routine application of Observation 4.1 and is accomplished by standard algorithms in linear algebra (reduction of the matrix to the canonical Jordan form).
\begin{itemize}
\item {\bf Step 1.} Compute the matrices $\alpha_i$, $\beta_i$
that constitute the matrix $M_{\rho_r}$ of the representation $\rho_r$.
\item {\bf Step 2.} Process the matrix of $M_{\rho_r}$ to derive 
the bar codes  ending up with a representation  
$\rho'_r$ whose all  $\alpha_i'$s  and $\beta_i'$s  are invertible matrices. 
\item {\bf Step 3} Compute the Jordan cells of $\rho_r$ 
from the representation $\rho_r'.$
\end{itemize}

\paragraph{Step 1.}
In Step 1 we begin with the incidence matrix
of the input simplicial complex $X$ equipped 
with the map $f:X\to {\mathbb S}^1.$ Let 
the angles $0\leq s_1 < s_2\cdots s_m \leq 2\pi$ be the critical values of $f$.
Choose a collection of regular angles $0 <t_1<t_2\cdots t_m <2\pi$ with 
$t_{i}< s_{i} <t_{i+1} < s_{i+1}$. Consider a canonical 
subdivision of $X$ into a cell complex so that $X_{[t_{i}, t_{i+1}]}$, 
and $X_{t_i}$ are subdivided into subcomplexes 
$R_i$ and $X_i$ as follows.
For any open simplex $\sigma$ we associate the open cells :
\begin{enumerate} 
 \item  $\sigma(i):=\sigma \cap X_{t_i}$ with $\dim (\sigma(i))= \dim \sigma-1$ if the intersection is nonempty
 \item  $\sigma \langle i\rangle:= \sigma\cap X_{(t_i, t_{i+1})}$ with
 $\dim \sigma \langle i \rangle= \dim \sigma$ if the intersection is nonempty.
\end{enumerate}

The cells of $X_i$ are exactly of the form $\sigma(i)$
and their incidences are given as $I(\sigma(i), \tau(i))= I(\sigma, \tau)$
where $I(\sigma,\tau)=0, +1, \mbox{ or } -1$ depending on whether $\tau$
is a coface of $\sigma$ and whether their orientations match or not.
The cells of $R_i$ consist of cells of $X_{i}$, $X_{i+1}$, and all
cells of the form $\sigma \langle i\rangle$. The incidences are given as 
$I(\sigma\langle i \rangle, \tau\langle i \rangle)= I(\sigma, \tau)$,  $I(\sigma(i),\sigma\langle i\rangle)=1$,
and  $I(\sigma(i+1), \sigma\langle i\rangle)=-1$. All other incidences are zero.
Assume that we are given a total order for the simplices of $X$ that is
compatible with $f$ and also the incidence relations.
This induces a total order for the cells in $X_i$ and $X_{i+1}$ and
also the cells in $R_i'=R_i\setminus X_{i}\sqcup X_{i+1}$ 
for any $1\leq i\leq m$ with $X_{m+1}:=X_1$. 
Impose a total order on $R_i$ by juxtaposing  the total orders of $X_{i}$, 
$X_{i+1}$, and $R_i'$ in this sequence. 
Clearly, the incidence matrix for $R_i$ 
can be derived from the incidence matrix of $X$. 

The incidence matrix of $A=X_i\sqcup X_{i+1}$ appears in the upper
left corner of the matrix for $R:=R_i$.
We obtain the matrices for the linear maps 
$\alpha_i:H_r(X_{t_{i}})\rightarrow H_r(X_{s_i})$ and
$\beta_i: H_r(X_{t_{i+1}})\rightarrow H_r(X_{s_{i}})$
by using the persistence algorithm~\cite{CEM06,ZC05}
on $R$ and $A$ as follows.
First, we run the persistence algorithm on the incidence 
matrix for $A$ to
compute a base of the homology group $H_r(A)$. We continue the procedure
by adding the columns and rows of the matrix for $R$ to obtain a base of
$H_r(R)$. It is straightforward to compute a matrix
representation of the inclusion induced linear map 
$H_r(A)\rightarrow H_r(R)$ with respect to the bases
computed by the persistence algorithm.

\paragraph{Step 2.}
Step 2 takes the matrix representation $M_{\rho_r}$ constructed out of 
matrices $\alpha_i, \beta_i$ computed in step 1,
and uses four elementary transformations $T_1(i),T_2(i),T_3(i)$, and $T_4(i)$ 
defined below to transform $M_{\rho_r}$ to  
$M_{\rho'_{r}}= T_{\cdots} (\cdots) M_{\rho_r},$ whose total 
number of rows and columns 
is strictly smaller than  that of $M_{\rho_r}$.
For convenience, let us write $\rho=\rho_r$ and $\rho'=\rho_r'$.
Each elementary transformation $T$ modifies the representation $\rho $ 
to the representation $\rho'$ while keeping indecomposable 
Jordan cell representations
unaffected but possibly changing the bar code representations. 
Some of these bar code representations remain the same, some are eliminated,
and some are shortened by one unit as 
described below. For each elementary transformation 
we record the changes to reconstruct the original bar codes. 
The elementary transformations are applied 
as long as the linear maps $\alpha_i$ or 
$\beta_i$ satisfy some injectivity and surjectivity property.
When no such transformation is applicable, the algorithm terminates with 
all $\alpha_i$ and $\beta_i$ being necessarily invertible matrices. 
At this point  the bar codes can be reconstructed reading  backwards the 
eliminations/modifications performed. The Jordan cells then can be 
obtained as detailed in Step 3. 

The elementary transformations modify the bar codes as follows: 
\begin{itemize}
\item $T_1(i)$ shortens  the bar codes$(i-1, k\}$ to $(i,k\}$  
if $i\geq 2$ and shortens the bar codes $(m, k\}, m<k$, 
to $(1, k-m\}$ if $i=1$.  
\item $T_2(i)$ shortens  the bar codes $\{l, i+km]$ to $\{l, i-1+km]$ 
for $k\geq 0.$ 
\item $T_3(i)$ shortens  the bar codes $[i, k\}$ to $[i+1,k\}$ 
for $i< m$ and to $[1, k-m\}$ if $i=m$. 
\item  $T_4(i)$ shortens  the bar codes $\{l, (i+1)+km )$ to $\{l, i+km)$ 
for $k\geq 0.$
\end{itemize}
It is understood that if an elementary transformation applied  to a bar code 
provides an interval which is not a bar code, then the bar code is eliminated.
Consequently  
$T_1(i)$ eliminates the bar codes $(i-1, i), (i-1, i]$\footnote{if $i=1$ 
eliminates the bar codes $(m,m+1)$ and $(m,m+1]$}, 
$T_2(i)$ eliminates the bar codes $[i,i], (i-1,i],$  
$T_3(i)$ eliminates the bar codes $[i,i+1), [i, i] $, and 
$T_4(i)$ eliminates the bar codes $(i, i+1), [i, i+1).$

To decide how many bar codes are eliminated  one 
uses Proposition \ref{P7} below.
Let $\#\{i,j\}_\rho$ denote the number of bar codes of type 
$\{i, j\}.$ 

\begin{proposition}\label {P7}
~
\begin{enumerate}
\item 
$\# (i,i+1)_\rho= \dim \ker\beta_i\cap\ker\alpha_{i+1}$
\item
 $\# [i,i]_\rho= \dim (V_{2i} / ((\beta_i (V_{2i+1}) + \alpha_i (V_{2i-1}))$
\item 
$\# (i,i+1]_\rho =
\dim(\beta_i(V_{2i+1}) +\alpha_i(\ker\beta_{i-1}))   -  \dim(\beta_i(V_{2i+1}))$
\item 
$\# [i, i+1)_\rho =
\dim(\alpha_i(V_{2i-1}) +\beta_i(\ker\alpha_{i+1}))-\dim(\alpha_i(V_{2i-1}))$
\end{enumerate}
\end{proposition}

The following diagrams define the elementary transformations 
and indicate the relation between the representation 
$\rho= \{V_i, \alpha_i, \beta_i\}$ and the representation 
$\rho' = \{V'_i, \alpha'_i, \beta'_i\}$ obtained after 
applying an elementary transformation.

\begin {itemize}
\item
Transformation $T_1(i)$:

\xymatrix{
&\cdots&V_{2i+1}\ar[l]_{\alpha_{i+1}}\ar [dr]_{\beta'_i}\ar[r]^{\beta_i}&V_{2i}\ar[d]              &V_{2i-1}\ar[l]_{\alpha_{i}}\ar[d]\ar[r]^{\beta_{i-1}}           &V_{2i-2}      &\cdots\ar[l]&\\
&           &                                                                                      &V'_{2i}                      &V'_{2i-1}\ar[l]^{\alpha'_i} \ar[ur]_{\beta'_{i-1}}        &                     &          &}

\vskip .2in 

\hskip 1in $V'_{2i-1}= V_{2i-1}/  \ker(\beta_{i-1}), \quad  V'_{2i}= V_{2i}/ \alpha_{i}(\ker(\beta_{i-1})),\quad V_k= V'_k, k\ne 2i, 2i-1$

\vskip .2in 

\item Transformation $T_2(i)$:

\xymatrix{
&\cdots&V_{2i+1}\ar[l]_{\alpha_{i+1}}\ar[dr]_{\beta'_{i}}\ar[r]^{\beta_i}&V_{2i}                      &V_{2i-1}\ar[l]_{\alpha_{i}}   \ar[r]^{\beta_{i-1}}           &V_{2i-2}      &\cdots\ar[l]&\\
&           &                               &V'_{2i}\ar[u]                      &V'_{2i-1}\ar[l]_{\alpha'_i} \ar[u]\ar[ur]_{\beta'_{i-1}}  &                     &          &}

\vskip .2in

\hskip 1in $V'_{2i}= \beta_i(V_{2i+1}),  \quad  V'_{2i-1}= \alpha^{-1}_{i}(\beta_i(V_{2i+1})), \quad V_k= V'_k, k\ne 2i-1, 2i$

\vskip .2in
\item Transformation $T_3(i)$:

\xymatrix{
&\cdots\ar[r]^{\beta_{i+1}} &V_{2i+2}  &V_{2i+1}\ar[l]_{\alpha_{i+1}}\ar[r] ^{\beta_i}             &V_{2i}                  &V_{2i-1}\ar[dl]^{\alpha'_i} \ar[l]_{\alpha_i}\ar[r]^{\beta_{i-1}}&\cdots&\\
&                    &                 &V'_{2i+1}\ar[ul]^{\alpha'_{i+1}}\ar[u]\ar[r]^{\beta'_i}   &V'_{2i}\ar[u]         &                                        &          &}
\vskip .2in 
\hskip 1in $V'_{2i}= \alpha_{i}(V_{2i-1}), \quad  V'_{2i+1}= \beta_i^{-1}(\alpha_i(V_{2i-1})), \quad   V_k= V'_k, k\ne 2i, 2i+1$
\vskip .2in 

\item Transformation $T_4(i)$:

\xymatrix{
&\cdots\ar[r]^{\beta_{i+1}}&V_{2i+2}&V_{2i+1}\ar[l]_{\alpha_{i+1}} \ar[d]\ar[r]^{\beta_i}&V_{2i}\ar[d]&V_{2i-1}\ar[l]_{\alpha_i} \ar[dl]^{\alpha'_i}\ar[r]^{\beta_{i-1}}&\cdots&\\
&                                           &                &V'_{2i+1}\ar[r]^{\beta'_i} \ar [ul]^{\alpha'_{i+1}}                &V'_{2i}         &                                        &          &}

\vskip .2in

\hskip 1in $V'_{2i+1}= V_{2i+1}/ \ker(\alpha_{i+1}), \quad  V'_{2i}= V_{2i}/ \beta_i(\ker(\alpha_{i+1})), \quad V_k= V'_k, k\ne 2i, 2i+1.$

\vskip .2in

\end{itemize}
The verification of the properties stated above and 
the proof of Proposition \ref{P7} are straightforward for  
indecomposable representations described in section \ref{S4} 
and  therefore for arbitrary representations.

As one can see from the diagrams above,
when  $\beta_{i-1}$  is injective, 
the representations $\rho$ and $\rho'$ are the same and we say 
that $T_1(i)$ is not applicable. 
Similarly, when $\beta_{i}$ is surjective, 
$T_2(i)$ is not applicable, when 
$\alpha_{i}$ 
is surjective, $T_3(i)$ is not applicable, and
when $\alpha_{i+1}$ is injective,  
$T_4(i)$ is not applicable. 
When  all $\alpha_i, \beta_i$ are invertible,
no elementary transformation is applicable 
and at this stage  the algorithm (Step 2) terminates.

To  explain how the algorithm works, it is convenient 
to  consider the following block matrices  
$B_{2i-1}$ and $B_{2i}$, $i=1,\cdots,m$, which
become the sub-matrices of
$M_{\rho_r}$ in (~\ref{main-matrix}) when the
entries $\beta_i$ are
replaced with $-\beta_i$. Let 
\begin{equation}\label{EEE1}
B_{2i-1}=
\begin{pmatrix}
\alpha_i& \beta_i\\
0& \alpha_{i+1}
\end{pmatrix},
\quad 
B_{2i}=
\begin{pmatrix}
\beta_i&0\\
\alpha_{i+1}& \beta_{i+1}
\end{pmatrix}
\end{equation}
 for $i= 1,2,\cdots (m-1)$ and 
 
 \begin{equation}\label{EEE2}
B_{2m-1}=
\begin{pmatrix}
\alpha_m& \beta_m\\
0& \alpha_{1}
\end{pmatrix},
\quad 
B_{2m}=
\begin{pmatrix}
\beta_m&0\\
\alpha_{1}& \beta_{1}
\end{pmatrix}.
\end{equation}
We modify $M_\rho$ by modifying successively each block $B_k.$ 
When $m>1$ 
the algorithm iterates over the blocks in multiple passes. In a single pass, 
it processes the blocks $B_1, B_2,\ldots,B_{2m}$ in this order.

When $B_{2(i-1)}= \begin{pmatrix}\beta_{i-1}&0\\ \alpha_i, &\beta_i\end{pmatrix}$ is processed then:
\begin {enumerate}
\item If $\beta_{i-1}$ is not injective, we apply  $T_1(i).$ 
This  boils down to changing the bases of $V_{2i-1}$ and  $V_{2i}$   
so that the matrix $B_{2(i-1)}$ becomes  
$$\left(\begin{array}{c|c||c}
\beta_{i-1,1}&0& 0\\ \hline\hline \alpha^1_{i,1}&\alpha^1_{i,2}&\beta^1_{i}\\ \hline \alpha^2_{i,1}&0& \beta^2_{i}\\\end{array}\right)$$
with $\begin{pmatrix}\beta_{i-1,1}&0\end{pmatrix}$ in column echelon form and $ \begin{pmatrix}\alpha^1_{i,2}\\0\end{pmatrix}$ in row echelon form. 

In this block matrix the first and third columns correspond to 
$V'_{2i-1}$ and $V_{2i+1}$ respectively, and the first and third rows to 
$V_{2(i-1)}$ and $V'_{2i}$ respectively. The second column and 
row become ``irrelevant" as a result of which
the modified block matrix $B_{2(i-1)}$ becomes 
 $\begin{pmatrix} \beta'_{i-1}&0\\ \alpha'_i&\beta'_i\end{pmatrix}= \begin{pmatrix}\beta_{i-1,1}&0\\ \alpha^2_{i,1}& \beta^2_{i}\end{pmatrix}.$

\item  If $\beta_{i}$ is not surjective, we apply $T_2(i).$ 
This  boils down to changing the bases of $V_{2i-1}$ and  $V_{2i}$   
so that the matrix $B_{2(i-1)}$ becomes  
$$\left(\begin{array}{c|c||c}
\beta_{i-1,1}&\beta_{i-1,2}& 0\\ \hline\hline \alpha^1_{i,1}&\alpha^1_{i,2}&\beta^1_{i}\\ \hline \alpha^2_{i,1}&0& 0\\ \end{array}\right)$$
with   $\begin{pmatrix} \beta^1_{i}\\0 \end{pmatrix}$ in row echelon form and 
$ \begin{pmatrix}\alpha^2_{i,1}& 0\end{pmatrix}$ in column echelon form. 

In this block matrix the second and third columns correspond to  
$V'_{2i-1}$ and $V_{2i+1}$ respectively, and the first and second rows to 
$V_{2(i-1)}$ and $V'_{2i}$ respectively.  
We make the first column and third row``irrelevant"
as a result of which
the modified block matrix $B_{2(i-1)}$ becomes 
 $\begin{pmatrix} \beta'_{i-1}&0\\ \alpha'_i&\beta'_i\end{pmatrix}=\begin{pmatrix}\beta_{i-1,2}&0\\ \alpha^1_{i,2}& \beta^1_{i} \end{pmatrix}.$
\end{enumerate}

When $B_{2i-1}$ is processed then:
\begin{enumerate}
\item[3.] If $\alpha_i$ is not surjective, 
we apply $T_3(i).$ This boils down to changing the bases 
of $V_{2i+1}$ and  $V_{2i}$   so that the matrix $B_{2i-1}$ becomes  
$$\left(\begin{array}{c||c|c}
\alpha^1_i&\beta^1_{i,1}& \beta^1_{i,2}\\ \hline0&\beta^2_{i,1}&0\\ \hline\hline 0&\alpha_{i+1,1}& \alpha_{i+1,2}\\ \end{array}\right)$$
with  
$\begin{pmatrix}\alpha^1_{i}\\0\end{pmatrix}$ 
in row echelon form and $ \begin{pmatrix}\beta^2_{i,1}&0\end{pmatrix}$ 
in column  echelon form. 

In this block matrix the first and third columns correspond to 
$V_{2i-1}$ and $V'_{2i+1}$ respectively, and the first and third rows to 
$V'_{2i}$ and $V_{2i+2}$ respectively.
We make the second column and second row ``irrelevant"
as a result of which
the modified block matrix $B_{2i-1}$ becomes 
 $\begin{pmatrix} \alpha'_{i}&\beta'_i\\0&\alpha'_{i+1}\end{pmatrix}= \begin{pmatrix}\alpha^1_{i}&\beta^1_{i,2}\\ 0 &\alpha_{i+1,2}\end{pmatrix}.$
\item[4.]  If $\alpha_{i+1}$ is not injective,
we apply $T_4(i).$ This boils down to changing the bases of 
$V_{2i+1}$ and  $V_{2i}$   so that the matrix $B_{2i-1}$ becomes  
$$\left(\begin{array}{c||c|c}
\alpha^1_i&\beta^1_{i,1}& \beta^1_{i,2}\\ \hline\alpha^2_i&\beta^2_{i,1}&0\\ \hline\hline 0&\alpha_{i+1,1}& 0\\ \end{array}\right)$$
with $\begin{pmatrix} \alpha_{i+1,1}&0 \end{pmatrix}$ 
in column echelon form and 
$ \begin{pmatrix}\beta^1_{i,2}\\ 0\end{pmatrix}$ in row echelon form. 

In this block matrix first and second columns correspond to 
$V_{2i-1}$ and $\ V'_{2i+1}$ respectively, and  second and third rows to
$V_{2i}'$ and $V_{2(i+1)}$ respectively. We make the third column and 
first row ``irrelevant" as a result of which
the modified block matrix $B_{2i-1}$ becomes 
 $\begin{pmatrix} \alpha'_{i}&\beta'_i\\0&\alpha'_{i+1}\end{pmatrix}= \begin{pmatrix}\alpha^2_{i}&\beta^2_{i,1}\\0& \alpha_{i+1,1} \end{pmatrix}.$
\end{enumerate}
Explicit formulae for $\alpha'$s  and $\beta'$s 
are given at the end of this section. 
At each pass the algorithm may eliminate or change 
bar codes, and if this happens, the matrix has less columns  or rows. 
If this does not happen, the algorithm terminates,
and indicates that there is no more bar code left. 
At termination, all $\alpha_i$ and $\beta_i$ become isomorphisms. 
The bar codes can be recovered by keeping track of all eliminations 
of the bar codes after each elementary transformation. A bar code
which is not eliminated in a pass gets shrunk by exactly two units,
during that pass, that is, a bar code $\{i,j\}$ shrinks to $\{i+1, j-1\}$ by exactly two
distinct elementary transformations. 
by elementary transformations.  For example if $m=5 $ the bar code $(1,5]$ during the pass became $(2,4]$ as result of 
applying $T_1(1)$ when inspecting $B_1$ and $T_2(5)$ when inspecting $B_9.$

When a bar code $[i,i]$ is eliminated, say, in
the $k$th pass, we know that it corresponds to a bar code $[i-k+1, i+k-1]$
in the original representation.
Similarly, other bar codes of type
$\{i,i+1\}$ eliminated at the $k$th pass correspond to the bar code
$\{i-k+1,i+k\}$. In both cases, the multiplicity of the bar codes
can be determined from the multiplicity of the eliminated bar codes thanks
to Proposition~\ref{P7}.

When $m=1$, the operations on above minors are not well defined. 
In this case we extend the quiver $G_2$ to $G_4$ ($m=2$) by adding fake
levels $t_2, s_2$ where $H_r(X_{t_2})=H_r(X_{s_2})=H_r(X_{s_1})$ and 
$\alpha_2,\beta_2$
are identities\footnote {Other easier methods can also be used in this case}.

A high level pseudocode for
the step~2 can be written as follows:
\vskip .2in

{\bf Algorithm} {\sc BarCode($M_{\rho}$)}\label {L}

\begin{quote}
\begin{enumerate}
\item[] Consider the block sub-matrices $B_1,\ldots,B_m$ of $M_{\rho}$;
\item[] Repeat
\begin{enumerate}
\item[] for $j:=1$ to $2m$ do
\begin{enumerate}
\item[1.] if $j=2i-1$ is odd 
\begin{enumerate}
\item if $\alpha_{i+1}$ is not injective, update
$B_{2i-1}:=T_4(i) (B_{2i-1})$.
\item if $\alpha_i$ is not surjective, update $B_{2i-1}:=T_3(i) (B_{2i-1})$.
\item delete any rows and columns rendered irrelevant.
\end{enumerate}
\item[2.] if $j=2i$ is even 
\begin{enumerate}
\item if $\beta_{i+1}$ is not surjective, update
$B_{2i}:=T_2(i)(B_{2i})$.
\item if $\beta_i$ is not injective, update $B_{2i}:=T_1(i)(B_{2i})$.
\item delete any rows and columns rendered irrelevant.
\end{enumerate}
\end{enumerate}
\item[]endfor
\end{enumerate}
\item[]until $M_\rho$ is not empty or
has not been updated.\\ 
Output $M_\rho$.
\end{enumerate}
\end{quote}

{\bf Example.} To illustrate how step 2 works,
we consider a representation given by 

\begin{equation}
\begin {aligned}
\alpha_1= \begin{pmatrix}1 & 1& 2\\
-3 &4& 2&\\ -2& 1&2\end{pmatrix} ;  
\alpha_2= \begin{pmatrix}1 & 0\\
0 &1\end{pmatrix};
\alpha_3= \begin{pmatrix}1&0& 0\\
0 &1&0\end{pmatrix};  \alpha_4= \begin{pmatrix}1 & 0\\
0 &1 \end{pmatrix}
%
\\
\beta_1= \begin{pmatrix} 1 & 0\\
0 &1 \\ 0& 0\end{pmatrix} ; 
\beta_2= \begin{pmatrix} 1 & 0& 0\\
0 &1 & 0\end{pmatrix};
\beta_3= \begin{pmatrix} 1 & 0\\
0 &1 \end{pmatrix}; \beta_4= \begin{pmatrix} 1&0& 0\\
0 &1 & 0\end{pmatrix} 
\end{aligned}
\end{equation}
The reader can notice that this is the representation $\rho_1$ for a simplified version of the example  provided in Fig 2 with the cylinder between the critical values $\theta_2$ and $\theta_3$ removed. 

\begin{itemize}
\item  Inspect $B_1$ and $B_2.$ No changes are necessary.
\item Inspect $B_3$.  Since $\alpha_3$ is not injective ,
one modifies the block by applying $T_4(2)$ which makes
both $\alpha_3$ and $\beta_2$ equal to $\begin{pmatrix} 1 & 0\\
0 &1 \end{pmatrix}$.
\item Inspect the blocks  $B_4, B_5, B_6, B_7.$ No changes are necessary. 
\item Inspect $B_8$. Since $\beta_4$ is not injective, one modifies 
the block by applying $T_1(1)$ which leads to 
$\alpha_1= \begin{pmatrix} -4 & 3\\
-3 &0 \end{pmatrix}$ and $\beta_1= \begin{pmatrix} -1 &1\\
-1 &0 \end{pmatrix}.$  
\end{itemize}
Indeed the block $B_8$ is given by  
$$ \left(\begin{array}{c||c}\beta_4&0\\\hline\hline \alpha_1&\beta_1\\ \end{array}\right)=  \left(\begin{array} {ccc||cc}
1 & 0& 0& 0& 0\\0 & 1& 0& 0& 0\\
\hline\hline 
1& 1& 2& 1& 0\\
-3 & 4& 2& 0&1\\-2 & 1& 2& 0& 0\\
\end{array}\right)$$

Since $\beta_4$ is already in column echelon form one only has to change the base of $V_2$ to bring the last column of $\alpha_1$ in row 
echelon form  which ends up with 

$$\left(\begin{array} {ccc||cc}
1 & 0& 0& 0& 0\\0 & 1& 0& 0& 0\\
\hline\hline 
1& 1& 2& 1& 0\\
-4 & 3& 0& -1&1\\-3 & 0& 0&-1& 0\\
\end{array}\right)$$
\vskip .1in
Therefore $\alpha'_1= \begin{pmatrix}-4& 3\\-3&0\end{pmatrix},$
$\beta'_1=\begin{pmatrix}1&0\\0&1\end{pmatrix},$
$\beta'_4=\begin{pmatrix}-1&1\\-1&0\end{pmatrix}.$
\vskip.1in

The algorithm stops  as all $\alpha_i'$s  and $\beta_i'$s are at this time invertible. 
 The last transformation  $T_1(1)$ has eliminated  only the bar code $(4,5]$,
 and the previous, which was the first  transformation, $T_4(2),$ has eliminated only  the bar code $(2, 3).$  This can be concluded from Proposition \ref{P7}.  In view of the properties of these 
two transformations, one concludes that these were the only two bar codes.

\paragraph{Step 3.}

At termination, all $\alpha_i$ and $\beta_i$ become isomorphisms because
otherwise one of the transformations would be applicable. The Jordan
cells can be recovered from the Jordan decomposition of the matrix
$$T= \beta^{-1}_{i-1}\cdot \alpha_{i-1}\cdot  \beta^{-1}_{i-2}\cdots \beta^{-1}_1\cdot \alpha_1\cdot \beta^{-1}_{m}\cdot \alpha_m\cdots    \beta^{-1}_{i+1}\cdot \alpha_{\i+1} \cdot \beta^{-1}_i\cdot \alpha_i \mbox{~~ for any } i.$$ Standard linear algebra routines permit the calculation of the Jordan cells   for familiar algebraic closed fields.
Note that if $\kappa$ is not algebraically closed,
Step 1 and Step 2 can still be performed and the matrix $T$ can be obtained. In this case it may not be possible to decompose the matrix $T$ in Jordan cells 
unless we consider the algebraic closure of $\kappa$. It is however possible
to decompose the matrix $T$ up to conjugacy as a sum of indecomposable 
invertible matrices while remaining in the class of matrices with coefficients
in the field $\kappa$.
This is the case for the field $\kappa= \mathbb Z_2$. 

\vskip.1in
 In the {\bf Example} above $T= \begin{pmatrix}3&1\\0&3\end{pmatrix}$  
provides the Jordan cell $(\lambda =3, k=2).$

\subsection{Implementation of $T_1(i), T_2(i), T_3(i)$ and $T_4(i).$}
\begin{enumerate}

\item[1.] $T_1 (i)$ acts on the block matrix 
$B_{2(i-1)}= \begin{pmatrix}\beta_{i-1}& 0\\ \alpha_i& \beta_i \end{pmatrix}$.
First we modify $B_{2(i-1)}$ to the block matrix 
$ \begin{pmatrix} \beta_{i-1,1}& 0&0\\ \alpha_{i,2}&\alpha_{i,2}& \beta_i \end{pmatrix}$ where
$\begin{pmatrix} \beta_{i-1,1}&0 \end{pmatrix}= \beta_{i-1}\cdot R(\beta_{i-1}) $ and   $\begin{pmatrix} \alpha_{i,1}&\alpha_{i,2} \end{pmatrix}= \alpha_i\cdot R(\beta_{i-1})$. Recall the definition of $R(\cdot)$ and $L(\cdot)$ given 
under notations in the introduction.
Then, one passes to the block matrix    
    $$\begin{pmatrix} \beta_{i-1,1}& 0&0\\ \alpha^1_{i,2}&\alpha^1_{i,2}& \beta^1_i\\ \alpha^2_{i,2}& 0& \beta^2_i\end{pmatrix} \mbox{ with } 
    \begin{pmatrix}\alpha^1_{i,2}\\ 0 \end{pmatrix} = L(\alpha_{i,2})\cdot \alpha_{i,2}, 
    \begin{pmatrix}\alpha^1_{i,1}\\ \alpha^2_{i,2} \end{pmatrix} = L(\alpha_{i,2}) \cdot \alpha_{i,1} \mbox{ and } 
     \begin{pmatrix}\beta^1_{i}\\ \beta^2_i \end{pmatrix} =L(\alpha_{i,2})\beta_{i}.$$
The modified block matrix is $ \begin{pmatrix}\beta_{i-1,1}&0\\ \alpha^2_{i,1}&\beta^2_i \end{pmatrix}.$

\item[2.] $T_2 (i)$   acts on the block matrix $B_{2(i-1)}= \begin{pmatrix}\beta_{i-1}& 0\\ \alpha_i& \beta_i \end{pmatrix}$.
First we modify $B_{2(i-1)}$ to the block matrix 
$ \begin{pmatrix} \beta_{i-1}& 0\\ \alpha^1_{i}& \beta^1_i\\ \alpha^2_i&0
 \end{pmatrix}$ where
 $\begin{pmatrix}\beta^1_{i}\\0\end{pmatrix}=L(\beta_i)\cdot \beta_i$ and 
 $\begin{pmatrix}\alpha^1_{i}\\ \alpha^2_i\end{pmatrix}=L(\beta_i)\cdot \alpha_i $
 
Then, one passes to the block matrix    
$$\begin{pmatrix} \beta_{i-1,1}& \beta_{i-1,2}&0\\ \alpha^1_{i,1}&\alpha^1_{i,2}& \beta^1_i\\ \alpha^2_{i,1}& 0&0\end{pmatrix} \mbox{ with } 
    \begin{pmatrix}\alpha^2_{i,1}\\ 0 \end{pmatrix} =   \alpha^2_{i,1}\cdot   R(\alpha^2_{i,1}),
    \begin{pmatrix}\alpha^1_{i,1}\\ \alpha^1_{i,2} \end{pmatrix} =  \alpha_{i,1}R(\alpha^2_{i,1}), \mbox{ and } 
    \begin{pmatrix}\beta_{i-1,1}\\ \beta_{i-1,2} \end{pmatrix} = \beta_{i-1}R(\alpha_{i,1}).$$
The modified block matrix is  $ \begin{pmatrix}\beta_{i-1,2}&0\\ \alpha^1_{i,2} &\beta^1_i \end{pmatrix}.$

\item[3.] $T_3 (i)$   acts on the block matrix $B_{2i-1}= \begin{pmatrix} \alpha_{i}& \beta_i\\ 0 & \alpha_{i+1} \end{pmatrix}$.
First we modify $B_{2i-1}$ to the block matrix 
$$ \begin{pmatrix} \alpha^1_i& \beta^1_i\\ 0 & \beta^2_i\\ 0&\alpha_{i+1}
 \end{pmatrix} \mbox{ where }
 \begin{pmatrix}\alpha^1_{i}\\0\end{pmatrix}=\alpha_i\cdot R(\alpha_i) 
\mbox{ and }
 \begin{pmatrix}\beta^1_{i}\\ \beta^2_i\end{pmatrix}=\beta_i\cdot R(\alpha_i).$$
Then, one passes to the block matrix    
    $$\begin{pmatrix} \alpha^1_i& \beta^1_{i,1}&\beta^1_{i,2}\\0 & \beta^2_{i,1}& 0\\ 0&\alpha_{i+1,1}&\alpha_{i+1,2}\end{pmatrix}
    \mbox{ with } 
    \begin{pmatrix}\beta^2_{i,1}&0\end{pmatrix} =  \beta^2_{i}\cdot   R( \beta^2_{i}),
    \begin{pmatrix}\beta^1_{i,1}& \beta^1_{i,2} \end{pmatrix} =  \beta^1_i\cdot   R( \beta^2_{i})$$
and $\begin{pmatrix}\alpha_{i+1,1}&\alpha_{i+1,2} \end{pmatrix} = \alpha_{i+1}\cdot   R( \beta^2_{i}).$
The modified block matrix is  $ \begin{pmatrix}\alpha^1_{i}&\beta^1_{i,2}\\ 0 &\alpha_{i+1,2}\end{pmatrix}.$

\item[4.] $T_4 (i)$   acts on the block matrix $B_{2i-1}= \begin{pmatrix} \alpha_{i}& \beta_i\\ 0 & \alpha_{i_1} \end{pmatrix}$.
First one modifies $B_{2i-1}$ to the block matrix 
$$ \begin{pmatrix} \alpha_i& \beta_{i,1}& \beta_{i,2} \\  0& \alpha_{i+1,1}& 0
 \end{pmatrix} \mbox{ where }
 \begin{pmatrix} \alpha_{i+1,1}& 0\end{pmatrix}=\alpha_{i+1}\cdot R(\alpha_{i+1}) \mbox{ and }
 \begin{pmatrix}\beta_{i,1} & \beta_{i,2} \end{pmatrix}=\beta_i\cdot R(\alpha_{i+1}).$$
 
Then, one passes to the block matrix    
    $$\begin{pmatrix} \alpha^1_i& \beta^1_{i,1}& \beta^1_{i,2} \\  \alpha^2_i& \beta^2_{i,1}& 0\\ 0 & \alpha^2_{i+1,1}& 0\end{pmatrix} \mbox{ with } 
    \begin{pmatrix}\beta^1_{i,2}\\ 0 \end{pmatrix} =   L(\beta_{i,2})\cdot   \beta_{i,2},
    \begin{pmatrix}\beta^1_{i,1}\\ \beta^2_{i,1} \end{pmatrix} =  L(\beta_{i,2})\cdot \beta_{i,1} \mbox{ and }
    \begin{pmatrix}\alpha^1_{i}\\ \alpha^2_{i} \end{pmatrix} = L(\beta_{i,2})\cdot \alpha_{i}.$$
The modified block matrix is  $ \begin{pmatrix}  \alpha^2_i& \beta^2_{i,1}\\ 0& \alpha_{i+1,1} \end{pmatrix}.$
\end{enumerate}

\subsection{Time complexity}
Let the input complex $X$ have $n$ simplices
in total on which the circle-valued map $f$ is defined
which has $m$ critical values. 

Then, step 1 takes $O(nd)$ time to detect all the critical values
where $d\leq n$ is the maximum degree of any vertex.
The critical values can be computed
by looking at the simplices adjacent to each of the vertices. To compute
the matrices $\alpha_i$ and $\beta_i$, we set up the matrices of
size $O(n)\times O(n)$ and run persistence on them. Using the algorithm
of~\cite{MMS11}, this can be achieved in $O(M(n))$ time
where $M(n)$ is the time complexity of multiplying two $n\times n$
matrices\footnote{We have $M(n)=O(n^{\omega})$
where $\omega < 2.376$
~\cite{CW90}.}.
Since we perform this operations for each of the critical levels 
and the spaces between them, we have $O(m M(n))$ total time complexity
for step 1.

In step 2, we process the matrix $M_{\rho_r}$ iteratively until all
BarCode representations are removed. In each pass except the
last one, we are guaranteed to shrink a bar code by at least
one unit. Therefore, the total number of passes is bounded
from above by the total length of all bar codes. Theorem~\ref{T1}
implies that a bar code cannot come back to the same level 
more than $\max_{s_i} \dim H_r(X_{s_i})$ times which can be
at most $O(n)$. Therefore, any bar code has a length of at most
$O(nm)$ giving a total length of $O(n^2m)$ over all bar codes.
Hence, the repeat loop in the algorithm {\sc BarCode} cannot
have more that $O(n^2m)$ iterations. In each iteration,
we reduce the block matrices each of which can be done with 
$O(M(n))$ matrix multiplication time
~\cite{Jean}. Since there are at most
$O(m)$ block matrices to be considered, we have $O(m M(n))$ time
per iteration giving a total of $O(n^2m^2 M(n))$ time for step 2.

Step 3 is performed on the resulting matrix from
step 2 which has $O(mn)\times O(mn)$ size.
This can again be performed by matrix multiplication
which takes $O(M(mn))$ time.

Therefore, the entire algorithm has time complexity of
$O(m^2n^2 M(n)+ M(mn))$.

\section{Conclusions} 
We have analyzed circle valued maps from the perspective of
topological persistence. We show that the notion of persistence
for such maps incorporate an invariant that is not encountered
in persistence studied erstwhile. Our results also shed lights
on computing  homology vector spaces and other topological invariants 
from bar codes and Jordan cells
(Theorems~\ref{T1} and~\ref{T2}). We have given an algorithm to
compute the bar codes and the Jordan cells; the algorithms can
also be adapted to compute zigzag persistence. 
In a subsequent work,  Burghelea and Haller 
have derived more subtle topological invariants like 
Novikov homology, monodromy~\cite {BH12},  Reidemeister torsion,
and others from bar codes and Jordan cells 
confirming  their mathematical relevance.
We have not treated in this paper the stability of the invariants;
see~\cite{BH12} for partial answer.

The standard persistence is related to Morse theory.
In a similar vein, the persistence for circle valued
map is related to Morse Novikov theory~\cite{Novikov}.
The work of Burghelea and Haller 
applies Morse Novikov theory to instantons and 
closed trajectories for vector field with 
Lyapunov closed one form~\cite{BH08}. 
The results in this paper 
will very likely provide additional insight on the 
dynamics of these vector fields and have implications 
in computational topology in particular and 
algebraic topology in general.

\section*{Acknowledgment} We acknowledge the support of the NSF grant
CCF-0915996 which made this research possible. We also thank all the referees
whose comments were helpful in improving the presentation of the paper.

\appendix
\section*{Appendix} \label {example}
 In this Appendix we explain the calculation of the $r$-invariants for 
the example depicted in Fig 2.
 The representation $\rho _0$ has vector spaces 
that are all one dimensional and maps $\alpha_i=
\beta_i$ that are all identity. 
Hence, there is no bar code, but one Jordan cell $\lambda=1, k=1.$ 

It is not hard to recognize from Fig 2  that  
the maps for the representation $\rho_1$  are given by:

\begin{equation*}
\begin {aligned}
\alpha_1= \begin{pmatrix}1 & 1&2\\
-3&4&2\\ -2&1&2\end{pmatrix} ; ~ \alpha_2=\begin{pmatrix} 1 & 0\\
0 &1\\ 0& 0\end{pmatrix} ;~ \alpha_3=\begin{pmatrix} 1 & 0& 0\\
0 &1 & 0\\ 0& 0&1\end{pmatrix};~ &\alpha_4= \begin{pmatrix}1 & 0\\
0 &1\end{pmatrix}\\
\alpha_5= \begin{pmatrix}1 & 0& 0\\
0 &1 & 0\end{pmatrix};~  &\alpha_6= \begin{pmatrix}1 & 0\\
0 &1 \end{pmatrix};\\
 \beta_1= \begin{pmatrix} 1 & 0\\
0 &1 \\ 0& 0\end{pmatrix} ;~ \beta_2= \begin{pmatrix} 1 & 0& 0\\
0 &1 & 0&\\ 0& 0& 1\end{pmatrix};~ \beta_3= \begin{pmatrix} 1 & 0\\
0 &1 \\ 0& 0\end{pmatrix};~ &\beta_4= \begin{pmatrix} 1 & 0& 0\\
0 &1 & 0\end{pmatrix}\\
\beta_5= \begin{pmatrix} 1 & 0\\
0 &1 \end{pmatrix};~ &\beta_6= \begin{pmatrix} 1 & 0& 0\\
0 &1 & 0 \end{pmatrix}.
\end{aligned}
\end{equation*}

We proceed with the step 2 of the algorithm.  
\begin {itemize}
\item inspect $B_1$ - no change for $\rho= \rho_1$;  inspect $B_2$- no change.
\item inspect $B_3$, - since $\alpha_2$ is not surjective apply $T_3(2).$ 
This changes $\alpha_2, \beta_2, \alpha_3$ into 
$\alpha'_2= \begin{pmatrix}1&0\\0&1\end{pmatrix},$ $\beta'_2= \begin{pmatrix}1&0\\0&1\end{pmatrix},$ $\alpha'_3= \begin{pmatrix}1&0\\0&1\\0&0\end{pmatrix}.$
Update and continue.
\item inspect $B_4$- no changes.
\item inspect $B_5$ - since $\alpha_3$ is not surjective, apply $T_3(3).$
This changes $\alpha_3$ and $\beta_3$ into $\alpha'_3= \begin{pmatrix}1&0\\0&1\end{pmatrix}$ and $\beta'_3= \begin{pmatrix}1&0\\0&1\end{pmatrix}.$
Update and continue.
\item inspect $B_6$- no changes.
\item inspect $B_7$ - since $\alpha_5$ is not injective, apply $T_4(4).$
This changes $\beta_4$ and $\alpha_5$ into $\alpha'_5= \begin{pmatrix}1&0\\0&1\end{pmatrix}$ and $\beta'_4= \begin{pmatrix}1&0\\0&1\end{pmatrix}.$
Update and continue.
\item inspect $B_8$ - no change; inspect $B_9$- no change; inspect $B_{10}$- no change; inspect $B_{11}$- no change.
\item inspect $B_{12}$ - since $\beta_6$ is not injective, apply $T_1(1).$
This changes $\beta_6$, $\alpha_1,$  $\beta_1$ to $\beta'_6= \begin{pmatrix}1&0\\0&1\end{pmatrix}$, $\alpha'_1= \begin{pmatrix}-4&3\\-3&0\end{pmatrix}$, and $\beta'_1= \begin{pmatrix}-1&1\\-1&0\end{pmatrix}.$
Update.
\end{itemize}
Since at this time all $\alpha_i'$s and $\beta_i'$s are invertible,
step 2 terminates.

\paragraph{Book keeping.} The last transformation $T_1(1)$  
has eliminated the bar code $(\theta_6, \theta_1+2\pi]$ (by Proposition \ref{P7}) and nothing else. This bar code was not the modification of any other bar code by the previous elementary transformations.  The previous transformation 
$T_4(4)$ has eliminated  the bar code $(\theta_4, \theta_5)$ and 
nothing else (by Proposition \ref{P7}). This bar code  was not the modification of any other bar code by the previous transformations.  The transformation $T_3(3)$has eliminated  the bar code  
$[\theta_3, \theta_3]$ (by Proposition \ref{P7}) which was the modification of $[\theta_2, \theta_3]$ by $T_3(2).$ These are all bar codes as listed 
in the table in section \ref{S3}. 
To calculate the Jordan cells we use step 3. We calculate the Jordan cells of $  \begin{pmatrix}-4&3\\-3&0\end{pmatrix}\cdot (\begin{pmatrix}-1&1\\-1&0\end{pmatrix})^{-1}$ which is  $(\lambda=3, k=2)$ as listed in the table in section \ref{S3}.

\vskip 1in

\end{document}

\begin{equation}
\begin {aligned}
\alpha_1= \begin{pmatrix}1 & 0& 0\\
0 &1 & 0&\\ 0& 0&1\end{pmatrix} ;  \alpha_2=\begin{pmatrix} 1 & 0\\
0 &1\\ 0& 0\end{pmatrix} ;\alpha_3=\begin{pmatrix} 1 & 0& 0\\
0 &1 & 0&\\ 0& 0&1\end{pmatrix}; \alpha_4= \begin{pmatrix}1 & 0\\
0 &1\end{pmatrix}\\
\alpha_5= \begin{vmatrix}1 & 0& 0\\
0 &1 & 0\end{vmatrix};  \alpha_6= \begin{vmatrix}1 & 0\\
0 &1 \end{vmatrix};  \alpha_7= \begin{vmatrix} 3 & 0& 2\\
6 &3 & 3&\\ 1& 1& 4\end{vmatrix} \\
\beta_1= \begin{vmatrix} 1 & 0\\
0 &1 \\ 0& 0\end{vmatrix} ; \beta_2= I\begin{vmatrix} 1 & 0& 0\\
0 &1 & 0&\\ 0& 0& 1\end{vmatrix}; \beta_3= \begin{vmatrix} 1 & 0\\
0 &1 \\ 0& 0\end{vmatrix}; \beta_4= \begin{vmatrix} 1 & 0& 0\\
0 &1 & 0\end{vmatrix}\\
\beta_5= \begin{vmatrix} 1 & 0\\
0 &1 \end{vmatrix}; \beta_6= ;\begin{vmatrix} 1 & 0& 0\\
0 &1 & 0\end{vmatrix} \beta_7=\begin{vmatrix} 1 & 0& 0\\
0 &1 & 0&\\ 0& 0& 1\end{vmatrix}
\end{aligned}
\end{equation}

The representation  $\rho_1$ has all $\alpha_i= \beta_i= Id_1,$  $ i=1\cdots 7$ and the representation $\rho_2$ is the trivial representation. There are 14 blocks defined by the formulae \ref{EEE1}

We inspect each of the blocks $B_1, B_2, \cdots B_{13}, B_{14},$  when needed make the change of he block, record the change and continue the inspection with the changes made.
No changes on $B_1$ and $B_2$ and change in $B_3$ in applying $T_3(B_3)$
Change $B_4$ applying $T_4(B_4)$ change again $B_4$ applying $T_1(B_4).$ 
no changes until  $B_7$ and change $B_7$ applying $T_2(B_7), $ no changes  until 
change $B_{12}$ by $T_1(B_{12})$ and change $B_{13}$ by $T_3(B_{13})$ and $B_{14}$ by   $T_4(B_{14})$ 

After the first pass  
we end up with a new representation $\rho_1'$ described by 

\begin{equation}
\begin {aligned} 
\alpha_1= \begin{vmatrix}1 & 0& 0\\
0 &1 & 0\end{vmatrix} ;  \alpha_2=\begin{vmatrix} 1 & 0\\
0 &1\end{vmatrix} ;\alpha_3=\begin{vmatrix} 1 & 0\\
0 &1 \end{vmatrix}; \alpha_4= \begin{vmatrix}1 & 0\\
0 &1\end{vmatrix}\\
\alpha_5= \begin{vmatrix}1 & 0\\
0 &1 \end{vmatrix};  \alpha_6= \begin{vmatrix}1 & 0\\
0 &1 \end{vmatrix};  \alpha_7= \begin{vmatrix} 3 & 0\\
6 &3 \end{vmatrix} \\
\beta_1= \begin{vmatrix} 1 & 0\\
0 &1 \end{vmatrix} ; \beta_2= \begin{vmatrix} 1 & 0\\
0 &1 \end{vmatrix}; \beta_3= \begin{vmatrix} 1 & 0\\
0 &1 \end{vmatrix}; \beta_4= \begin{vmatrix} 1 & 0\\
0 &1 \end{vmatrix}\\
\beta_5= \begin{vmatrix} 1 & 0\\
0 &1 \end{vmatrix}; \beta_6= ;\begin{vmatrix} 1 & 0\\
0 &1 \end{vmatrix} \beta_7=\begin{vmatrix} 1 & 0& 0\\
0 &1 & 0&\\ 0& 0& 1\end{vmatrix}
\end{aligned}
\end{equation}

\end{document}

\paragraph{Elementary transformations.}
\label{jordan}
In our algorithm, we used four elementary 
transformations $T_1, T_2, T_3$, and $T_4$
that modify the matrix form of a quiver representation.
We elaborate on these four transformations here. First, we 
interpret these transformations in terms of the quiver representations,
and then give their implementations in terms of matrices.

Each transformation takes an index $i$ and a representation 
$$\rho=\{V_j| 1\le j\le 2m , \, \alpha_s: V_{2s-1}\to V_{2s}, \, 
\beta_s: V_{2s+1}\to V_{2s}| 1\le s\le m\}$$ 
and produces a new representation 
$$\rho'=\{V'_j| 1\le j \le 2m ,\,  \alpha'_s: V'_{2s-1}\to V'_{2s}, \, \beta'_s: V'_{2s+1}\to V'_{2s}| 1\le s\le m\}$$ 
as follows:
\begin{enumerate}
\item 
If $\rho'= T_1(\rho,i)$  then  $V'_{2i-1}= V_{2i-1}/ \ker (\beta_{i-1})$,  $V'_{2i}=V_{2i}/\alpha_{i}(\ker(\beta_{i-1})),$  $V'_{j}= V_j$  for $j\ne \{2i-1, 2i\}$ with $\alpha'_s, \beta'_s$ being induced from $\alpha_s, \beta_s$ for
$s\in [1,m]$.
\item 
If  $\rho'= T_2(\rho,i)$ then  $V'_{2i+1}= V_{2i+1}/ \ker (\alpha_{i+1})$,  $V'_{2i}=V_{2i}/\beta_{i}(\ker\alpha_{i+1}),$  $V'_{j}= V_j$  for $j\ne \{2i+1, 2i\}$ with $\alpha'_s, \beta'_s$ being induced from $\alpha_s, \beta_s$ for 
$s\in[1,m]$.
\item 
If $\rho'= T_3(\rho,i)$  then $V'_{2i}= \alpha_i(V_{2i-1})$, $V'_{2i+1}=\beta_i^{-1}(\alpha_i(V_{2i-1})),$  $V'_{j}= V_j$ for $j\ne \{2i, 2i+1\}$ with $\alpha'_s, \beta'_s$ being the restrictions  of $\alpha_s, \beta_s$ for $s\in[1,m]$.
\item 
If $\rho'= T_4(\rho,i)$ then $V'_{2i}= \beta_i(V_{2i+1})$, $V'_{2i-1}=\alpha_i^{-1}(\beta_i(V_{2i+1}))$,  $V'_{j}= V_j$ for $j\ne \{2i, 2i-1\}$ with $\alpha'_s, \beta'_s$ being the restrictions  of $\alpha_s, \beta_s$ for $s\in [1,m]$.
\end{enumerate}

The following diagrams \footnote{in the diagrams $V_0= V_{2m} $ 
and $\beta_0= \beta_m$} indicate the constructions described above. 
The indices increase
from right to left to signify that the vector spaces are laid counterclockwise
with increasing indices around a quiver.\\

Transformation $T_1(\rho,i)$:

\xymatrix{
&\cdots&V_{2i+1}\ar[l]_{\alpha_{i+1}}\ar [dr]_{\beta'_i}\ar[r]^{\beta_i}&V_{2i}\ar[d]              &V_{2i-1}\ar[l]_{\alpha_{i}}\ar[d]\ar[r]^{\beta_{i-1}}           &V_{2i-2}      &\cdots\ar[l]&\\
&           &                                                                                      &V'_{2i}                      &V'_{2i-1}\ar[l]^{\alpha'_i} \ar[ur]_{\beta'_{i-1}}        &                     &          &}

\vskip .2in 

\hskip 1in $V'_{2i-1}= V_{2i-1}/  \ker(\beta_{i-1}) \quad  V'_{2i}= V_{2i}/ \alpha_{i}(\ker(\beta_{i-1}))$

\vskip .2in 

Transformation $T_2(\rho,i)$:

\xymatrix{
&\cdots\ar[r]^{\beta_{i+1}}&V_{2i+2}&V_{2i+1}\ar[l]_{\alpha_{i+1}} \ar[d]\ar[r]^{\beta_i}&V_{2i}\ar[d]&V_{2i-1}\ar[l]_{\alpha_i} \ar[dl]^{\alpha'_i}\ar[r]^{\beta_i}&\cdots&\\
&                                           &                &V'_{2i+1}\ar[r]^{\beta'_i} \ar [ul]^{\alpha'_{i+1}}                &V'_{2i}         &                                        &          &}

\vskip .2in

\hskip 1in $V'_{2i+1}= V_{2i+1}/ \ker(\alpha_{i+1}), \quad  V'_{2i}= V_{2i}/ \beta_i(\ker(\alpha_{i+1}))$

\vskip .2in

Transformation $T_3(\rho,i)$:

\xymatrix{
&\cdots\ar[r]^{\beta_{i+1}} &V_{2i+2}  &V_{2i+1}\ar[l]_{\alpha_{i+1}}\ar[r] ^{\beta_i}             &V_{2i}                  &V_{2i-1}\ar[dl]^{\alpha'_i} \ar[l]_{\alpha_i}\ar[r]^{\beta_i}&\cdots&\\
&                    &                 &V'_{2i+1}\ar[ul]^{\alpha'_{i+1}}\ar[u]\ar[r]^{\beta'_i}   &V'_{2i}\ar[u]         &                                        &          &}
\vskip .2in 
\hskip 1in $V'_{2i}= \alpha_{i}(V_{2i-1}) \quad  V'_{2i+1}= \beta_i^{-1}(\alpha_i(V_{2i-1}))$
\vskip .2in 

Transformation $T_4(\rho,i)$:

\xymatrix{
&\cdots&V_{2i+1}\ar[l]_{\alpha_{i+1}}\ar[dr]_{\beta'_{i}}\ar[r]^{\beta_i}&V_{2i}                      &V_{2i-1}\ar[l]_{\alpha_{i}}   \ar[r]^{\beta_{i-1}}           &V_{2i-2}      &\cdots\ar[l]&\\
&           &                               &V'_{2i}\ar[u]                      &V'_{2i-1}\ar[l]_{\alpha'_i} \ar[u]\ar[ur]_{\beta'_{i-1}}  &                     &          &}

\vskip .2in

\hskip 1in $V'_{2i}= \beta_i(V_{2i+1})  \quad  V'_{2i-1}= \alpha^{-1}_{i}(\beta_i(V_{2i+1}))$

\vskip .1in

The following observations follow straightforwardly from the
definitions.

$T_1(\rho,i)$  eliminates all bar codes of the form 
$(i-1,i)$ and $(i-1,i],$
shrinks each bar code of the form $(i-1,k]$ and $(i-1,k),$ $k\geq (i+2),$ 
into bar codes $(i,k]$ and $(i,k)$ respectively
with the convention that $(i-1,k\}$ is $(m, m+k\}$ when $i=1$, and
leaves all  other bar codes and Jordan cells unchanged. 

$T_2(\rho,i)$ eliminates all bar codes of the 
form $(i,i+1)$ and $[i,i+1),$
shrinks each bar code of the form $[l,i+1)$ and $(l,i+1),$  $l\leq i-1,$ 
into bar codes $[l,i)$ and $(l,i)$ respectively,
and leaves any other bar codes and Jordan cells unchanged.

Type $T_3(\rho,i)$ 
eliminates all bar codes of the form $[i,i]$ and $[i,i+1),$
shrinks each bar code of the forms $[i, k)$ and $[i,k]$, $k\geq i+1,$ into 
the bar codes $[i+1,k)$ and $[i+1, k]$ respectively,
and leaves all  other  type of bar codes and Jordan cells unchanged.

Type $T_4(\rho,i)$ 
eliminates all bar codes of the form $[i,i]$ and $(i-1,i],$ 
shrinks each bar code of the forms $(l, i] $ and $[l,i],$ $l\leq i-1,$ 
into the bar codes $(l, i-1]$ and $[l,i-1]$ respectively 
with the convention that $\{l,0\}$ is identified to $\{l+m, m\}$, 
and leaves all other type of bar codes and Jordan cell unchanged.
%

Let $\#\{i,j\}_\rho$ denote the number of bar codes of type 
$\{i, j\}$ \footnote{ we use $"\{"$ resp. for $"\}"$ as a common notation for $"["$ and $"( "$ resp. $"\}"$ and $")"$.}  for a representation $\rho.$
We have the following:

\begin{proposition}\label {P7}
~
\begin{enumerate}
\item 
$\# (i,i+1)_\rho= \dim \ker\beta_i\cap\ker\alpha_{i+1}$
\item
 $\# [i,i]_\rho= \dim (V_{2i} / ((\beta_i (V_{2i+1}) + \alpha_i (V_{2i-1}))$
\item 
$\# (i,i+1]_\rho =
\dim(\beta_i(V_{2i+1}) +\alpha_i(\ker\beta_{i-1}))   -  \dim(\beta_i(V_{2i+1}))$
\item 
$\# [i, i+1)_\rho =
\dim(\alpha_i(V_{2i-1}) +\beta_i(\ker\alpha_{i+1}))-\dim(\alpha_i(V_{2i-1}))$
\end{enumerate}
\end{proposition}
which can be easily derived from the matrices $\alpha_i$ and $\beta_i$. 

\paragraph{Implementing the elementary transformations.}
Elementary transformations are implemented with algorithms
known to produce the so-called {\it echelon form} of a matrix.
We briefly describe the concept of echelon form for completeness.
Let $A$ be any $m\times n$ matrix.
For any column $j$ in $A$, denote by $r(j)\in \{1,\cdots ,m\}$ 
the first index so that  $a_{r(j),j}\ne 0$, i.e., $a_{s,j}=0$ if $s< r(j).$
We say that $A$ is in  column- echelon form if:
\begin{enumerate}
\item $a_{r(j),j}=1$ 
\item $r(j) < r(j+1)$
\item $a_{r(j),s} =0$ for $s<j$
\end{enumerate} and  in row echelon form
 if the transpose $A^\ast$ is in column echelon form. 
The following is well-known in linear algebra.

\begin{obs}
Given a matrix $A$,
\begin{enumerate}
\item there exists an
$m\times m$ matrix $L(A)$ so that $L(A)\cdot A$ is in row echelon form;
\item there exists an $n\times n$ matrix $R(A)$ 
so that $A\cdot R(A)$ is in column echelon form;
\item the rank and the dimension of the kernel of $A$ can be 
read off from its echelon form.
\end{enumerate}
\end{obs}

There are well known algorithms to obtain the matrices $L(A)$ and $R(A)$
to produce the echelon forms.
 
Now we describe how we take advantage of the echelon forms to
implement the elementary transformations. 
In transformations $T_2$ and $T_3$, we replace the matrix 
\begin{equation*}
\begin{vmatrix}
\alpha_i& -\beta_i\\
0& \alpha_{i+1}
\end{vmatrix}
\rm{ first\  by }
\begin{vmatrix}
\bar\alpha_i& -\bar\beta_i\\
0& \bar\alpha_{i+1}
\end{vmatrix}
\rm{and \ then \  by}
\begin{vmatrix}
\bar{\bar\alpha}_i& -\bar{\bar\beta}_i\\
0& \bar{\bar\alpha}_{i+1}
\end{vmatrix}
\end{equation*}

while in transformations $T_1$ and $T_4$ we replace 
 
\begin{equation*}
\begin{vmatrix}
 -\beta_i &0\\
 \alpha_{i+1}& -\beta_{i+1}
\end{vmatrix}
\rm{ first\  by }
\begin{vmatrix}
-\bar\beta_i &0\\
 \bar\alpha_{i+1}& -\bar\beta_{i+1}
\end{vmatrix}
\rm{and \ then \  by}
\begin{vmatrix}
-\bar{\bar\beta}_i &0\\
 {\bar{\bar\alpha}}_{i+1}& -\bar{\bar\beta}_{i+1}
\end{vmatrix}.
\end{equation*}

Because of that we will also write 
$$ T_1(B_{2i}), \ T_4(B_{2i}), \ T_2(B_{2i-1}), \ T_3(B_{2i-1})$$
for 
$$T_1(\rho, i+1),\ T_4(\rho, i+1),\ T_2(\rho, i),\ T_3(\rho, i).$$

The blocks (submatrices) in the second and the third matrices represent the 
same linear transformations 
as the ones in the first one but with respect to different bases of the 
same vector spaces $V_i$s.
The columns and rows which became zero by passing to the second and 
to third forms are referred to as ``irrelevant" and  
are  removed 
without damaging the Jordan cells but  implying elimination or  changes 
in the bar codes as specified earlier.\\

\noindent
{\bf Transformation $T_2(\rho,i)=T_2(B_{2i-1})=T_2
\begin{vmatrix}
\alpha_i& -\beta_i\\
0& \alpha_{i+1}
\end{vmatrix}$}:\\

The matrices below summarize the modifications carried out for 
this transformation.
\vskip 0.1in 
\[
\begin{array}{ccccc}

\left( \begin{array}{c|c}
\alpha_i & \beta_i \\ \hline
0 & \alpha_{i+1}
\end{array}
\right)
&
\Rightarrow
&
\left( \begin{array}{c|c|c}
\begin{array}{c}
\bar{\alpha}_i=\alpha_i\\ \hline 
0
\end{array} 
&
\begin{array}{c}
\bar{\beta}_{i,1}\\ \hline
\bar\alpha_{i+1,1} 
\end{array}
&
\begin{array}{c}
\bar\beta_{i,2}\\ \hline 
0
\end{array} 
\end{array}
\right)
&
\Rightarrow
&
\left( \begin{array}{c|c|c}
\begin{array}{c}
\begin{array}{c}
{\bar{\bar{\alpha}}}_{i,1}\\ \hline 
{\bar{\bar{\alpha}}}_{i,2}
\end{array}\\ \hline 
0
\end{array}

&

\begin{array}{c}
\begin{array}{c}
\bar{\bar{\beta}}_{i,1,1}\\ \hline
{\bar{\bar{\beta}}}_{i,1,2}
\end{array}\\ \hline
{\bar{\bar\alpha}}_{i+1,1} 
\end{array}
&
\begin{array}{c}
\begin{array}{c}
{\bar{\bar\beta}}_{i,2}\\ \hline 
0
\end{array}\\ \hline
0
\end{array}
\end{array}
\right)
\end{array}
\] 

\noindent
First take 
\begin{equation*}
\begin{aligned}
\bar\alpha_i=& \alpha_i\\
\bar\beta_i= &\beta_i\cdot R(\alpha_{i+1})=(\bar\beta_{i,1},\ \bar\beta_{i, 2})\\
\bar\alpha_{i+1}=& \alpha_{i+1}\cdot R(\alpha_{i+1})=(\bar\alpha_{i+1,1},\ 0)
\end{aligned}
\end{equation*}
and then take
\begin{equation*}
\begin{aligned}
\bar{\bar\alpha}_{i}=& L(\bar\beta_{i,2})\cdot\bar\alpha_{i}=
\begin{vmatrix}
\bar{\bar\alpha}_{i,1}\\
\bar{\bar\alpha}_{i,2} 
\end{vmatrix}\\
\bar{\bar\beta}_{i,1}= &L(\bar\beta_{i,2})\cdot \bar\beta_{i,1} =
\begin{vmatrix}
\bar{\bar\beta}_{i,1,1}\\
\bar{\bar\beta}_{i,1,2}
\end{vmatrix}
\\
\bar{\bar\beta}_{i,2}= &L(\bar\beta_{i,2})\cdot \bar\beta_{i,2} =  
\begin{vmatrix}
\bar{\bar\beta}_{i,2,1}\\
0
\end{vmatrix}\\
\bar{\bar\alpha}_{i+1}= &\bar\alpha_{i+1}\Rightarrow \bar{\bar\alpha}_{i+1,1}=\bar\alpha_{i+1,1}
\end{aligned}
\end{equation*}
The ``irrelevant columns" correspond to the last elements of the modified 
base of $V_{2i+1}$ so that the zero columns of $\bar\alpha_{i+1}$ vanish.
The ``irrelevant rows" correspond to  the last elements of the modified 
base of $V_{2i}$ so that the zero rows of $L(\bar\beta_{i,2})\cdot \bar\beta_{i,2}$ 
vanish. 

The minor after the elimination of the irrelevant rows and columns becomes 
$$\boxed{
\begin{vmatrix}
\alpha'_i \ \ &\beta'_i\\
0\ \ &\alpha'_{i+1}
\end{vmatrix}:=
\begin{vmatrix}
\bar{\bar\alpha}_{i,2} \ \ &\bar{\bar\beta}_{i,1, 2}\\
0\ \ &\bar{\bar\alpha}_{i+1,1}
\end{vmatrix}}.
$$


\noindent
{\bf Transformation $T_3(\rho, i)= T_3(B_{2i-1})=T_3
\begin{vmatrix}
\alpha_i& -\beta_i\\
0& \alpha_{i+1}
\end{vmatrix}$}:\\

\noindent
First take 
\begin{equation*}
\begin{aligned}
\bar\alpha_i=& L(\alpha_i)\cdot  \alpha_i=
\begin{vmatrix}
\bar\alpha_{i,1}\\
0
\end{vmatrix} \\
\bar\beta_i= &L(\alpha_i)\cdot \beta_i=
\begin{vmatrix}
\bar\beta_{i,1}\\
\bar\beta_{i,2}
\end{vmatrix}\\
\bar\alpha_{i+1}= &\alpha_{i+1}
\end{aligned}
\end{equation*}
and then take
\begin{equation*}
\begin{aligned}
\bar{\bar\alpha}_i= &\bar\alpha_i \Rightarrow \bar{\bar\alpha}_{i,1}=\bar\alpha_{i,1}\\
\bar{\bar\beta}_{i,1}= &\bar\beta_{i,1}\cdot R(\bar\beta_{i,2}) =(\bar{\bar\beta}_{i,1,1},\  \bar{\bar\beta}_{i,1,2})\\
\bar{\bar\beta}_{i,2}= &\bar\beta_{i,2} \cdot R(\bar\beta_{i,2})= (\bar{\bar\beta}_{i,2,1},\  0)  \\
\bar{\bar\alpha}_{i+1}=& \bar\alpha_{i+1}\cdot R(\bar\beta_{i,2})=(\bar{\bar\alpha}_{i+1,1},\  \bar{\bar\alpha}_{i+1,2}) 
\end{aligned}
\end{equation*}

The minor after the elimination of the irrelevant rows and columns becomes 
$$\boxed{
\begin{vmatrix}
\alpha'_{i} &  \beta'_i\\
0 & \alpha'_{i+1}
\end{vmatrix}:=
\begin{vmatrix}
\bar{\bar\alpha}_{i,1} &  \bar{\bar\beta}_{i,1,2}\\
0 & \bar{\bar\alpha}_{i+1,2}
\end{vmatrix}
}$$


\noindent
{\bf Transformation $T_4(\rho, i+1)= T_4(B_{2i})=T_4
\begin{vmatrix}
 -\beta_i &0\\
 \alpha_{i+1}& -\beta_{i+1}
\end{vmatrix}$}:\\

\noindent
First take 
\begin{equation*}
\begin{aligned}
\bar\beta_i= &\beta_i\\
\bar\alpha_{i+1}= &L(\beta_{i+1})\cdot \alpha_{i+1}=
\begin{vmatrix}
\bar\alpha_{i+1,1}\\
\bar\alpha_{i+1,2}
\end{vmatrix}\\
\bar\beta_{i+1}=&L(\beta_{i+1})\cdot \beta_{i+1}=
\begin{vmatrix}
\bar\beta_{i+1,1}\\
0
\end{vmatrix}
\end{aligned}
\end{equation*}
and then take
\begin{equation*}
\begin{aligned}
\bar{\bar\beta}_i= &\bar\beta_{i}\cdot R(\bar\alpha_{i+1,2})=(\bar{\bar\beta}_{i,1}\ \bar{\bar\beta}_{i,2}) \\
\bar{\bar\alpha}_{i+1,1}= &\bar\alpha_{i+1,1}\cdot R(\bar\alpha_{i+1,2}) =(\bar{\bar\alpha}_{i+1,1,1}\  \bar{\bar\alpha}_{i+1,1,2})\\
\bar{\bar\alpha}_{i+1,2}= &\bar\alpha_{i+1,2}\cdot R(\bar\alpha_{i+1,2}) =(\bar{\bar\alpha}_{i+1,2,1}\  0)\\
\bar{\bar\beta}_{i+1}=& \bar\beta_{i+1} \Rightarrow 
\bar{\bar\beta}_{i+1,1}=\bar\beta_{i+1,1}\\
\end{aligned}
\end{equation*}

The minor after the elimination of the irrelevant rows and columns is 
$$\boxed{
\begin{vmatrix}
\beta'_{i} \ \ & 0\\
\alpha'_{i+1}\ \ &\beta'_{i+1}
\end{vmatrix}:=
\begin{vmatrix}
\bar{\bar\beta}_{i,1} \ \ & 0\\
\bar{\bar\alpha}_{i+1,1,2}\ \ &\bar{\bar\beta}_{i+1,2}
\end{vmatrix}
}.$$
%
\vskip .3in
\noindent
{\bf Transformation $T_1(\rho, i+1)= T_1(B_{2i})=T_1
\begin{vmatrix}
 -\beta_i &0\\
 \alpha_{i+1}& -\beta_{i+1}
\end{vmatrix}$}:\\

\noindent
First take 
\begin{equation*}
\begin{aligned}
\bar\beta_{i}=& \beta_{i}\cdot R(\beta_{i})=(\bar\beta_{i,1},\ 0)
\\
\bar\alpha_{i+1}=& \alpha_{i+1}\cdot R(\beta_i)=(\bar\alpha_{i+1,1},\ \bar\alpha_{i+1, 2})
\\
\bar\beta_{i+1}=&\beta_{i+1}
\end{aligned}
\end{equation*}
and then take

\begin{equation*}
\begin{aligned}
\bar{\bar\beta}_{i}= &\bar\beta_{i} \Rightarrow \bar{\bar\beta}_{i,1}=\bar\beta_{i,1} \\
\bar{\bar\alpha}_{i+1,2}= &L(\bar\alpha_{i+1,2})\cdot \bar\alpha_{i+1,2} =
\begin{vmatrix}
\bar{\bar\alpha}_{i+1,2,1}\\
0
\end{vmatrix}
\\
\bar{\bar\alpha}_{i+1,1}= &L(\bar\alpha_{i+1,2})\cdot \bar\alpha_{i+1,1} =  
\begin{vmatrix}
\bar{\bar\alpha}_{i+1,1,1}\\
\bar{\bar\alpha}_{i+1,1,2}
\end{vmatrix}\\
\bar{\bar\beta}_{i+1}=& L(\bar\alpha_{i+1,2})\cdot\bar\beta_{i+1}=
\begin{vmatrix}
\bar{\bar\beta}_{i+1,1}\\
\bar{\bar\beta}_{i+1, 2} 
\end{vmatrix}
\end{aligned}
\end{equation*}
The minor after the elimination of the irrelevant rows and columns becomes 
$$\boxed{
\begin{vmatrix}
\beta'_{i}  \ \ \ 0\\
\alpha'_{i+1} \   \beta'_{i+1}
\end{vmatrix}:=
\begin{vmatrix}
\bar{\bar\beta}_{i,1}  \ \ \ 0\\
\bar{\bar\alpha}_{i+1,1,2} \   \bar{\bar\beta}_{i+1,2}
\end{vmatrix}
}.$$

\begin{equation*}
\begin {aligned}
\alpha_1= \begin{vmatrix}1 & 1&2\\
-3&4&2\\ -2&1&2\end{vmatrix} ;  \alpha_2=\begin{vmatrix} 1 & 0\\
0 &1\\ 0& 0\end{vmatrix} ;\alpha_3=\begin{vmatrix} 1 & 0& 0\\
0 &1 & 0\\ 0& 0&1\end{vmatrix}; &\alpha_4= \begin{vmatrix}1 & 0\\
0 &1\end{vmatrix}\\
\alpha_5= \begin{vmatrix}1 & 0& 0\\
0 &1 & 0\end{vmatrix};  &\alpha_6= \begin{vmatrix}1 & 0\\
0 &1 \end{vmatrix};\\
 \beta_1= \begin{vmatrix} 1 & 0\\
0 &1 \\ 0& 0\end{vmatrix} ; \beta_2= \begin{vmatrix} 1 & 0& 0\\
0 &1 & 0&\\ 0& 0& 1\end{vmatrix}; \beta_3= \begin{vmatrix} 1 & 0\\
0 &1 \\ 0& 0\end{vmatrix}; &\beta_4= \begin{vmatrix} 1 & 0& 0\\
0 &1 & 0\end{vmatrix}\\
\beta_5= \begin{vmatrix} 1 & 0\\
0 &1 \end{vmatrix}; &\beta_6= ;\begin{vmatrix} 1 & 0& 0\\
0 &1 & 0 \end{vmatrix}
\end{aligned}
\end{equation*}

Here is what the Step 2 of the algorithm does:
\begin {itemize}
\item inspect $B_1$ - no change;  inspect $B_2$- no change;
\item inspect $B_3$, - since $\alpha_2$ is not surjective apply $T_3(2)$ 
This changes $\alpha_2, \beta_2, \alpha_3$ into 
$\alpha'_2= \begin{vmatrix}1&0\\0&1\end{vmatrix},$ $\beta'_2= \begin{vmatrix}1&0\\0&1\end{vmatrix},$ $\alpha'_3= \begin{vmatrix}1&0\\0&1\\0&0\end{vmatrix}.$
Update $M_\rho$ and continue.
\item inspect $B_4$- no changes;
\item inspect $B_5$ - since $\alpha_3$ is not surjective apply $T_3(3)$
This changes $\alpha_3$ and $\beta_3$ into $\alpha'_3= \begin{vmatrix}1&0\\0&1\end{vmatrix}$ and $\beta'_3= \begin{vmatrix}1&0\\0&1\end{vmatrix}.$
Update and continue.
\item inspect $B_6$- no changes;
\item inspect $B_7$= since $\alpha_5$ is not injective  apply $T_2(4).$
This changes $\beta_4$ and $\alpha_5$ into $\alpha'_5= \begin{vmatrix}1&0\\0&1\end{vmatrix}$ and $\beta'_4= \begin{vmatrix}1&0\\0&1\end{vmatrix}.$
Update and continue.
\item inspect $B_8$ - no change; inspect $B_9$- no change; inspect $B_{10}$- no change; inspect $B_{11}$- no change;
\item inspect$B_12$- since $\beta_6$ is not injective apply $T_1(1).$
This changes $\beta_6$ $\alpha_1,$  $\beta_1$ into $\beta'_6= \begin{vmatrix}1&0\\0&1\end{vmatrix}$ $\alpha'_1= \begin{vmatrix}-4&3\\-3&0\end{vmatrix}$ and $\beta'_1= \begin{vmatrix}-1&1\\-1&0\end{vmatrix}.$
Update.
\end{itemize}
Since at this time all $\alpha_i'$s and $\beta_i'$s are invertible the algorithm for Step 2 stops.

{\bf Book keeping.} The last transformation $T_1(1)$ has eliminated only  $(\theta_6, \theta_1+\2\pi]$ (by Proposition \ref{P7}) which was not changed by the previous transformations.The previous transformation 
$T_2(4)$ has eliminated only the bar code $(\theta_4, \theta_5)$ (by Proposition \ref{P7}) which was not modified by the previous transformation. The transformation transformation $T_3(3)$has eliminated only the bar code  
$[\theta_3, \theta_3]$ (by Proposition \ref{P7}) which was the modification of $[\theta_2, \theta_3]$ by $T_3(2).$ These are all bar codes.
To calculate the Jordan cells one uses Step 3 and calculates the Jordan cells of $( \begin{vmatrix}-1&1\\-1&0\end{vmatrix})^{-1}\cdot \begin{vmatrix}-4&3\\-3&0\end{vmatrix}$ which is only one, $(\lambda=3, k=2).$

\section{Conclusions}
We have analyzed circle valued maps from the perspective of
topological persistence. We show that the notion of persistence
for such maps incorporate an invariant that is not encountered
in persistence studied erstwhile. Our results also shed lights
on computing  homology vector spaces and other topological invariants 
from bar codes and Jordan cells
(Theorems~\ref{T1} and~\ref{T2}). We have given an algorithm to
compute the bar codes and the Jordan cells; the algorithms can be 
also adapted to compute zigzag persistence. 
In subsequent work  Burghelea and Haller 
have derived more subtle topological invariants like 
Novikov homology, monodromy~\cite {BH12},  Reidemeister torsion,
and others from bar codes and Jordan cells 
confirming  their mathematical relevance.
We have not treated in this paper the stability of the invariants. The issue will be addressed in a forthcoming paper, see~\cite{BH12} for partial answer.

The standard persistence is related to Morse theory.
In a similar vein, the persistence for circle valued
map is related to Morse Novikov theory~\cite{Novikov}.
The work of Burghelea and Haller 
applies Morse Novikov theory to instantons and 
closed trajectories for vector field with 
Lyapunov closed one form~\cite{BH08}. 
The results in this paper 
will very likely provide additional insight on the 
dynamics of these vector fields~\cite{BH08} and have implications 
in computational topology in particular and 
algebraic topology in general.

\section*{Acknowledgment} We acknowledge the support of the NSF grant
CCF-0915996 which made the collaboration between the two authors on this
research possible.

\section*{Appendix} \label {example}
After applying Step1 of the algorithm to the example described in section \ref{S3} one obtains the representation $\rho_0, \rho_1$ of the graph $G_{14}.$
The representation $\rho_0$ has all vector spaces one dimensional and all $\alpha_i=
\beta_i$ the identity, hence leads only to a Jordan cell $(1;1)$ in dimension zero.  The representation $\rho_1,$ as can be recognized  in Fig 2  is given by :

\begin{equation*}
\begin {aligned}
\alpha_1= \begin{vmatrix}1 & 1&2\\
-3&4&2\\ -2&1&2\end{vmatrix} ;  \alpha_2=\begin{vmatrix} 1 & 0\\
0 &1\\ 0& 0\end{vmatrix} ;\alpha_3=\begin{vmatrix} 1 & 0& 0\\
0 &1 & 0\\ 0& 0&1\end{vmatrix}; &\alpha_4= \begin{vmatrix}1 & 0\\
0 &1\end{vmatrix}\\
\alpha_5= \begin{vmatrix}1 & 0& 0\\
0 &1 & 0\end{vmatrix};  &\alpha_6= \begin{vmatrix}1 & 0\\
0 &1 \end{vmatrix};\\
 \beta_1= \begin{vmatrix} 1 & 0\\
0 &1 \\ 0& 0\end{vmatrix} ; \beta_2= \begin{vmatrix} 1 & 0& 0\\
0 &1 & 0&\\ 0& 0& 1\end{vmatrix}; \beta_3= \begin{vmatrix} 1 & 0\\
0 &1 \\ 0& 0\end{vmatrix}; &\beta_4= \begin{vmatrix} 1 & 0& 0\\
0 &1 & 0\end{vmatrix}\\
\beta_5= \begin{vmatrix} 1 & 0\\
0 &1 \end{vmatrix}; &\beta_6= ;\begin{vmatrix} 1 & 0& 0\\
0 &1 & 0 \end{vmatrix}
\end{aligned}
\end{equation*}

Here is what the Step 2 of the algorithm does:
\begin {itemize}
\item inspect $B_1$ - no change;  inspect $B_2$- no change;
\item inspect $B_3$, - since $\alpha_2$ is not surjective apply $T_3(2)$ 
This changes $\alpha_2, \beta_2, \alpha_3$ into 
$\alpha'_2= \begin{vmatrix}1&0\\0&1\end{vmatrix},$ $\beta'_2= \begin{vmatrix}1&0\\0&1\end{vmatrix},$ $\alpha'_3= \begin{vmatrix}1&0\\0&1\\0&0\end{vmatrix}.$
Update $M_\rho$ and continue.
\item inspect $B_4$- no changes;
\item inspect $B_5$ - since $\alpha_3$ is not surjective apply $T_3(3)$
This changes $\alpha_3$ and $\beta_3$ into $\alpha'_3= \begin{vmatrix}1&0\\0&1\end{vmatrix}$ and $\beta'_3= \begin{vmatrix}1&0\\0&1\end{vmatrix}.$
Update and continue.
\item inspect $B_6$- no changes;
\item inspect $B_7$= since $\alpha_5$ is not injective  apply $T_2(4).$
This changes $\beta_4$ and $\alpha_5$ into $\alpha'_5= \begin{vmatrix}1&0\\0&1\end{vmatrix}$ and $\beta'_4= \begin{vmatrix}1&0\\0&1\end{vmatrix}.$
Update and continue.
\item inspect $B_8$ - no change; inspect $B_9$- no change; inspect $B_{10}$- no change; inspect $B_{11}$- no change;
\item inspect$B_12$- since $\beta_6$ is not injective apply $T_1(1).$
This changes $\beta_6$ $\alpha_1,$  $\beta_1$ into $\beta'_6= \begin{vmatrix}1&0\\0&1\end{vmatrix}$ $\alpha'_1= \begin{vmatrix}-4&3\\-3&0\end{vmatrix}$ and $\beta'_1= \begin{vmatrix}-1&1\\-1&0\end{vmatrix}.$
Update.
\end{itemize}
Since at this time all $\alpha_i'$s and $\beta_i'$s are invertible the algorithm for Step 2 stops.

{\bf Book keeping.} The last transformation $T_1(1)$ has eliminated only  $(\theta_6, \theta_1+\2\pi]$ (by Proposition \ref{P7}) which was not changed by the previous transformations.The previous transformation 
$T_2(4)$ has eliminated only the bar code $(\theta_4, \theta_5)$ (by Proposition \ref{P7}) which was not modified by the previous transformation. The transformation transformation $T_3(3)$has eliminated only the bar code  
$[\theta_3, \theta_3]$ (by Proposition \ref{P7}) which was the modification of $[\theta_2, \theta_3]$ by $T_3(2).$ These are all bar codes.
To calculate the Jordan cells one uses Step 3 and calculates the Jordan cells of $( \begin{vmatrix}-1&1\\-1&0\end{vmatrix})^{-1}\cdot \begin{vmatrix}-4&3\\-3&0\end{vmatrix}$ which is only one, $(\lambda=3, k=2).$


\end{document}